\documentclass[leqno]{article}
\usepackage{latexsym,amsmath,amsfonts,mathrsfs,amssymb,amsthm}
\usepackage{yhmath}
\usepackage[usenames]{color}
\usepackage{graphicx}
\def\R{\mathbb R}

\def\Z{\mathbb Z}

\def\a{\alpha}
\def\e{\epsilon}
\def\d{\delta}
\def\Y{\mathbb Y}
\def\T{\mathbb T}

\def\be{\begin{equation}}
\def\ee{\end{equation}}
\def\bs{\backslash}
\def\qed{\hfill$\Box$\bigskip}
\def\nd{\noindent Proof. }
\def\F{{\cal F} }
\def\K{{\cal K} }
\numberwithin{equation}{section}
\newtheorem*{varthm}{Theorem 5.1}
\newtheorem{lem}[equation]{Lemma}
\newtheorem{pro}[equation]{Proposition}
\newtheorem{defn}[equation]{Definition}
\newtheorem{thm}[equation]{Theorem}
\newtheorem{cor}[equation]{Corollary}
\newtheorem{rem}[equation]{Remark}
\begin{document}
\bigskip

\centerline{\Large \textbf{Almgren and topological minimality for the set $Y\times Y$}}

\bigskip

\centerline{\large Xiangyu Liang}

\vskip 1cm

\centerline {\large\textbf{Abstract.}}

In this paper we discuss various minimality properties for the orthogonal product of two 1-dimensional $\Y$ sets, and some related problems. This is motivated by an attempt to give the classification of singularities for 2-dimensional Almgren-minimal sets in $\R^4$.
\bigskip

\textbf{AMS classification.} 28A75, 49Q10, 49Q20, 49K99

\bigskip

\textbf{Key words.} Minimal sets, Classification of singularities, Calibration, Topological decomposition, Hausdorff measure.

\setcounter{section}{-1}

\section{Introduction and preliminaries}

The main purpose of this paper is to prove the Almgren minimality of the orthogonal product of two one dimensional $\Y$ sets. This is part of the classification of singularities for 2-dimensional Almgren-minimal sets in $\R^4$. 
%
%In \cite{LM94} the authors exploited the beautiful method of paired calibration to proving minimality for sets of codimension 1, in particular the 1-dimensional $\Y$ set. While we tried to use the product of  the paired calibration to prove the product of two minimal set, a serious difficulty is that the codimension of the product is larger than one: separation condition no longer exists, and deformations can have more complicated behavior so that the calibration can easily fail to calibrate well. To overcome this issue we use the homology group on the set (rather than on the complementary, as the separation condition) to decompose the set properly. The calculation of the calibration becomes the estimation of the mass norm of non-simple 2-vectors.

The notion of Almgren minimality was introduced by Almgren to modernize Plateau's problem, which aims at understanding physical objects, such as soap films, that minimize the area while spanning a given boundary.  The study of regularity and existence for these sets is one of the centers of interest in geometric measure theory.

Recall that  a most simple version (without singularities) of Plateau's problem aims at finding a surface which minimizes area among all the surfaces having a given curve as boundary. See works of Besicovitch, Federer, Fleming, De Giorgi, Douglas, Reifenberg, etc. for example. Lots of other notions of minimality have also been introduced to modernize Plateau's problem, such as mass minimizing or size minimizing currents (see \cite{DPHa} for their definitions), varifolds (cf.\cite{Al66}). In this article, we will mainly use the notion of minimal sets introduced by F.Almgren \cite{Al76}, in a general setting of sets, and which gives a very good description of the behavior of soap films. 

Soap films that interested Plateau are 2-dimensional objects, but a general definition of $d-$dimensional minimal sets in an open set $U\subset \R^n$ is not more complicated. 

Intuitively, a $d-$dimensional minimal set $E$ in an open set $U\subset\R^n$ is a closed set $E$ whose $d-$dimensional Hausdorff measure could not be decreased by any local Lipschitz deformation. The more precise definition is the following:

\begin{defn}[Almgren competitor (Al competitor for short)] Let $E$ be a relatively closed set in an open subset $U$ of $\R^n$ and $d\le n-1$ be an integer. An Almgren competitor for $E$ is a relatively closed set $F\subset U$ that can be written as $F=\varphi_1(E)$, where $\varphi_t:U\to U$ is a family of continuous mappings such that 
\be \varphi_0(x)=x\mbox{ for }x\in U;\ee
\be\mbox{ the mapping }(t,x)\to\varphi_t(x)\mbox{ of }[0,1]\times U\mbox{ to }U\mbox{ is continuous;}\ee
\be\varphi_1\mbox{ is Lipschitz,}\ee
  and if we set $W_t=\{x\in U\ ;\ \varphi_t(x)\ne x\}$ and $\widehat W=\bigcup_{t\in[0.1]}[W_t\cup\varphi_t(W_t)]$,
then
\be \widehat W\mbox{ is relatively compact in }U.\ee
 
Such a $\varphi_1$ is called a deformation in $U$, and $F$ is also called a deformation of $E$ in $U$.
\end{defn}

\begin{defn}[Almgren minimal sets]
Let $0<d<n$ be integers, $U$ an open set of $\R^n$. A relatively closed set $E$ in $U$ is said to be Almgren minimal of dimension $d$ in $U$ if 
\be H^d(E\cap B)<\infty\mbox{ for every compact ball }B\subset U,\ee
and
\be H^d(E\bs F)\le H^d(F\bs E)\ee
for all Al competitors $F$ for $E$.
\end{defn}

\bigskip

The point of view here is very different from those of minimal surfaces and mass minimizing currents under certain boundary conditions, which are more often used. Comparing to the large number of results in the theory of mass minimizing currents, or classical minimal surfaces, in our case, very few results of regularity and existence are known. However, Plateau's problem is more like the study of size minimizing currents, for which likewise very few results are known (see \cite{DP09} for certain existence results). One can prove that the support of a size minimizing current is automatically an Almgren minimizer, so that all the regularity results listed below are also true for supports of size minimizing currents.

\medskip

The first regularity results for minimal sets have been given by Frederick Almgren \cite{Al76} (rectifiability, Ahlfors regularity in arbitrary dimension), then generalized by Guy David and Stephen Semmes \cite{DS00} (uniform rectifiability, big pieces of Lipschitz graphs), Guy David \cite{GD03} (minimality of the limit of a sequence of minimizers). 

Since minimal sets are rectifiable and Ahlfors regular, they admit a tangent plane at almost every point. But our main interest is to study those points where there is no tangent plane, i.e. singular points.

\medskip

A first finer description of the interior regularity for minimal sets is due to Jean Taylor, who gave in \cite {Ta} an essential regularity theorem for 2-dimensional minimal sets in 3-dimensional ambient spaces: if $E$ is a minimal set of dimension 2 
in an open set of $\R^3$ , then  every point $x$ of $E$ has a neighborhood where $E$ is equivalent (modulo a negligible set) through a $C^1$ diffeomorphism to a minimal cone (that is, a minimal set which is also a cone). 

In \cite{DJT}, Guy David generalized Jean Taylor's theorem to 2-dimensional minimal sets in $\R^n$, but with a local bi-H\"older equivalence, that is, every point $x$ of $E$ has 
a neighborhood where $E$ is equivalent through a bi-H\"older diffeomorphism to a minimal cone $C$ (but the minimal cone might not be unique). 

In addition, in \cite{DEpi}, David also proved that, if this minimal cone $C$ satisfies a ``full-length" condition, we will have the $C^1$ equivalence (called $C^1$ regularity). In particular, the tangent cone of $E$ at the point $x\in E$ exists and is a minimal cone, and the blow-up limit of $E$ at $x$ is unique; if the full-length condition fails, we might lose the $C^1$ regularity.

Thus, the study of singular points is transformed into the classification of singularities, i.e., into looking for a list of minimal cones. Besides, getting such a list would also help deciding locally what kind (i.e. $C^1$ or bi-H\"older) of equivalence with a minimal cone can we get.

\bigskip

In $\R^3$, the list of 2-dimensional minimal cones has been given by several mathematicians a century ago. (See for example \cite{La} or \cite{He}). They are, modulo isometry: a plane, a $\Y$ set (the union of 3 half planes that meet along a straight line where they make angles of 120 degrees), and a $\T$ set (the cone over the 1-skeleton of a regular tetrahedron centered at the origin). See the pictures below.

\centerline{\includegraphics[width=0.2\textwidth]{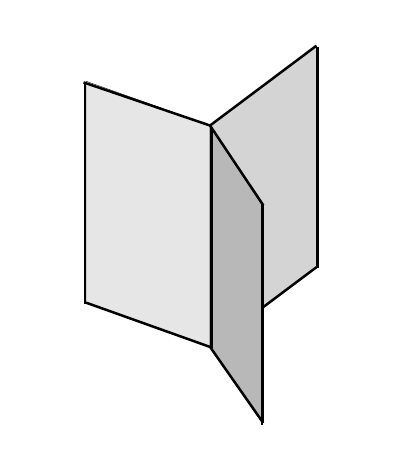} \hskip 2cm\includegraphics[width=0.25\textwidth]{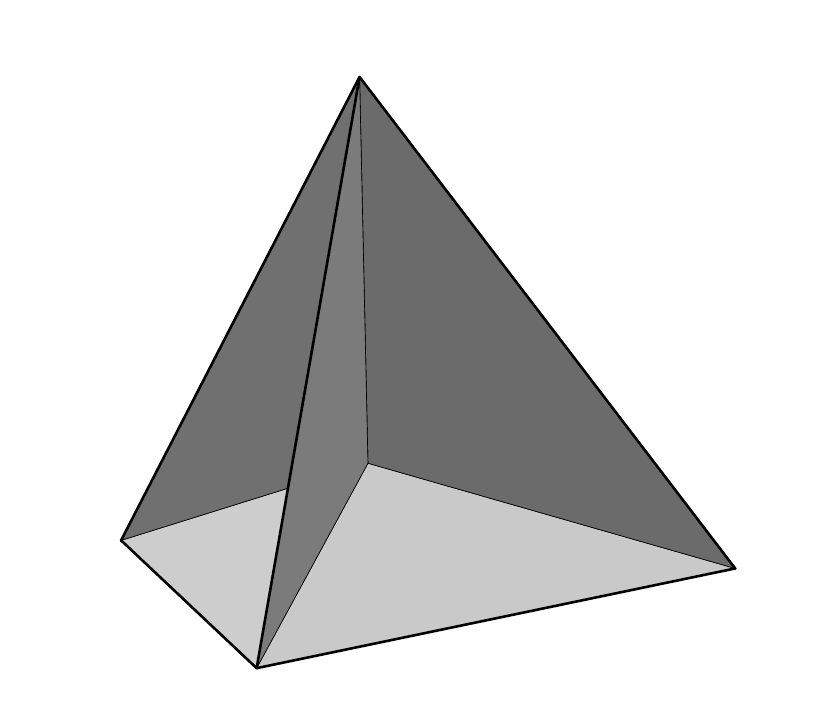}}

In higher dimensions, even in dimension 4, the list of minimal cones is still very far from clear. Except for the three minimal cones that already exist in $\R^3$, the only 2-dimensional minimal cones that were known before this paper is the union of two almost orthogonal planes (cf.\cite{2p} Thm 1.2). But this continuous one-parameter family of minimal cones gives already an interesting phenomenon in $\R^4$ that does not occur in $\R^3$----in $\R^3$, no small perturbation of any minimal cones ever preserves the minimality, and moreover, each minimal cone admits a different topology from the others. 

The next natural and only known candidate for a 2-dimensional minimal cone in $\R^4$ is probably the set $Y\times Y$, that is, the orthogonal product of two 1-dimensional $\Y$ sets (the union of three half lines with a common endpoint $x$, that meet at $x$ at $120^\circ$ angles).  The main theorem of this paper is the following (see Section 5):

\begin{varthm}The set $Y\times Y$ is an Almgren minimal cone in $\R^4$.
\end{varthm}

Notice that the 1-dimensional $\Y$ sets are 1-dimensional minimal cones in $\R^n$. In fact it seems very natural to think that the product of any two Almgren minimal sets is still Almgren-minimal. However in general, we do not know how to prove it, even though our Theorem \ref{al} affirms that this is true for the particular case of $Y\times Y$. We even do not know whether the product of an Almgren minimal set with $\R$ is minimal or not. This inconvenience might be due to the fact that the Almgren minimality is the weakest among all other related notions of minimality. For instance, if we take the notion of topological minimality (cf.\cite{topo}), at least the product of a topological minimal set $E$ with $\R$ is still topological minimal.
%
%This is true for the product of two 1-dimensional $\Y$ set, due to our result. But in general, we do not know how to prove that the product of two minimal sets is still minimal. We even do not know whether the product of an Almgren-minimal set $E$ with $\R$ is still Almgren minimal. This inconvenience might be due to the fact that the Almgren minimality is the weakest among all other related notions of minimality. For instance, if we take the notion of topological minimality (cf.\cite{topo}), at least the product of a topological minimal set $E$ with $\R$ is still topological minimal.

Conversely, we can prove that if the orthogonal product set $E=E_1\times E_2$ is Almgren minimal, then each $E_i,i=1,2$ is minimal. We will deal with this in the last section.

\bigskip

The proof of Theorem \ref{al} uses the particular topological structure of the one dimensional $\Y$ sets. Recall that Brakke \cite{Br91},  Lawlor and Morgan \cite{LM94} have introduced the beautiful technique of paired calibrations for proving the minimality of various sets of codimension 1, based on some separation condition. Note that the separation condition is preserved by deformations, it helps to decompose any Almgren competitor for $\Y$ into pieces, each of which can be calibrated well by the paired calibrations.

%
%While we tried to use the product of  the paired calibration to prove the product of two minimal set, a serious difficulty is that the codimension of the product is larger than one: separation condition no longer exists, and deformations can have more complicated behavior so that the calibration can easily fail to calibrate well. To overcome this issue we use the homology group on the set (rather than on the complementary, as the separation condition) to decompose the set properly. The calculation of the calibration becomes the estimation of the mass norm of non-simple 2-vectors.

For the set $Y\times Y$, we still want to use the product of the calibrations, but we face a serious difficulty and new ideas are required:  the set $Y\times Y$ is of codimension 2, where the separation condition no longer exists, and deformations can have more complicated behavior so that the product of the calibrations can easily fail to calibrate well if we still decompose an Almgren competitor according to the associated deformation. In fact in her thesis \cite{XY10}, when the author tried to use the product of the paired calibrations for the one dimensional $\Y$ set to calibrate also the image of any deformation of the product of two one-dimensional $\Y$ sets, she only managed to prove that $Y\times Y$ is just minimal among the class of all its injective deformations. When a deformation is not injective, there might be some unexpected intersections of images of different parts of $Y\times Y$, where the calibration can easily fail to calibrate well. 

In this paper we deal with a deformation of $Y\times Y$ differently, by decomposing it according to some new topological condition on its homology group. The point is that topological properties are somehow intrinsic for a set, while the deformation is relatively unstable, in the sense that one set can have many different parameterizations. Note also that in codimension 1, the separation condition is a condition imposed on the complements of sets; but here our topological condition is imposed on the set itself.

For this purpose we will first define a larger class $\F$ of competitors (see Section 1) for a set $E$, which are sets that satisfy this new topological condition. And then a big part (Section 1-4) of the paper is devoted to showing that $Y\times Y$ is minimal in this class of competitors $\F$, by decomposing any set in $\F$ according to this topological condition. We then show (in Section 5) that the class of Almgren competitors, as well as another class of topological competitors, are contained in the big class $\F$. As a result, $Y\times Y$ minimizes naturally the Hausdorff measure in these two smaller classes of competitors, and thus is both Almgren minimal and topological minimal.

But unlike the union of two almost orthogonal planes in $\R^4$, for which the minimality is stable under small perturbations, we will prove that our set $Y\times Y$ is an isolated type of singularity. That is, a small perturbation of $Y\times Y$ could not give a minimal cone. We know even more, that is, no other 2-dimensional minimal cone in $\R^4$ could admit the same topology as $Y\times Y$. 
%
%It is also very natural to ask whether the product of any two Almgren minimal sets is Almgren-minimal. This is true for the product of two 1-dimensional $\Y$ set, due to our result. But in general we do not know how to prove it. We even do not know whether the product of an Almgren-minimal set $E$ with $\R$ is still Almgren minimal. This inconvenience might be due to the fact that the Almgren minimality is the weakest among all other related notions of minimality. For instance, if we take the notion of topological minimality (cf.\cite{topo}), at least the product of a topological minimal set $E$ with $\R$ is still topological minimal.

\bigskip

The plan for the rest of this article is the following. 

In Section 1 we introduce the class $\F$ of fundamental competitors that preserve some topological condition, and try to establish the equivalence of the minimality in this class, and the minimality of the subclass $\F_R$ of regular fundamental competitors, that are decomposable.

In Section 2 we give the decomposition of regular competitors, and prove some intersection and projection properties for these decompositions. The aim of the decomposition is to avoid those ugly intersections of pieces, so that the sum of the calibrations could not be too large.

In Section 3 we give the calibration, and prove the minimality of $Y\times Y$ in the class $\F$.

In Section 4 we prove the maximum of the calibration (Lemma 3.24).

In Section 5 we establish the Almgren and topological minimality of $Y\times Y$.

We discuss perturbations of $Y\times Y$ in Section 6, and prove that no other 2-dimensional Almgren minimal cone in $\R^4$ could admit the same topology as $Y\times Y$.

In Section 7 we prove that for general sets the Almgren minimality of $E_1\times E_2$ yields the Almgren minimality for $E_i,i=1,2$.

\noindent\textbf{Some useful notation}

$[a,b]$ is the line segment with end points $a$ and $b$;

$[a,b)$ is the half line with initial point $a$ and passing through $b$;

$B(x,r)$ is the open ball with radius $r$ and centered on $x$;

$\overline B(x,r)$ is the closed ball with radius $r$ and center $x$;

$\overrightarrow{ab}$ is the vector $b-a$;

$H^d$ is the Hausdorff measure of dimension $d$.

\section{Regular fundamental competitors in $D$}

Write $\R^4=\R^2_1\times \R^2_2$. For $i=1,2$, $Y_i\subset \R^2_i$ is a 1-dimensional $\Y$ set centered at the origin of $\R^2_i$. Denote by $Y\times Y$ the set $Y_1\times Y_2\subset \R^4$.

For any $r>0$, and $i=1,2$, denote by $B_i(x,r)\subset\R^2_i$ the open ball centered at $x\in \R^2_i$ of radius $r$. In particular we denote by $B_i,i=1,2$ the unit ball in $\R^2_i$. Denote by $C_i\subset \R^2_i$ the closed convex hull of $Y_i\cap B_i$, then $C_i$ is a regular triangle inscribed in $B_i$. Set $D=C_1\times C_2$.

Denote by $a_j, j=1,2,3$, the three points of intersection of $Y_1$ with $\partial B_1$, and by $b_j,j=1,2,3$, the three points of intersection of $Y_2$ with $\partial B_2$. 

Denote by $S_1$ the union of segments $[a_1,o_1]\cup [o_1,a_2]$, $S_2=[a_2,o_1]\cup [o_1,a_3]$, and $S_3=[a_3,o_1]\cup [o_1,a_1]$; similarly set $R_1=[b_1,o_2]\cup [o_2,b_2]$, $R_2=[b_2,o_2]\cup [o_2,b_3]$, $R_3=[b_3,o_2]\cup [o_2,b_1]$, where $o_i,i=1,2$ denotes the origin of $\R^2_i$. Then $S_i\subset Y_1\subset\R_1^2$, $1\le i\le 3$, and  $R_j\subset Y_2\subset\R_2^2$,  $1\le j\le 3$. (See Figure 1 below).

 \centerline{\includegraphics[width=0.7\textwidth]{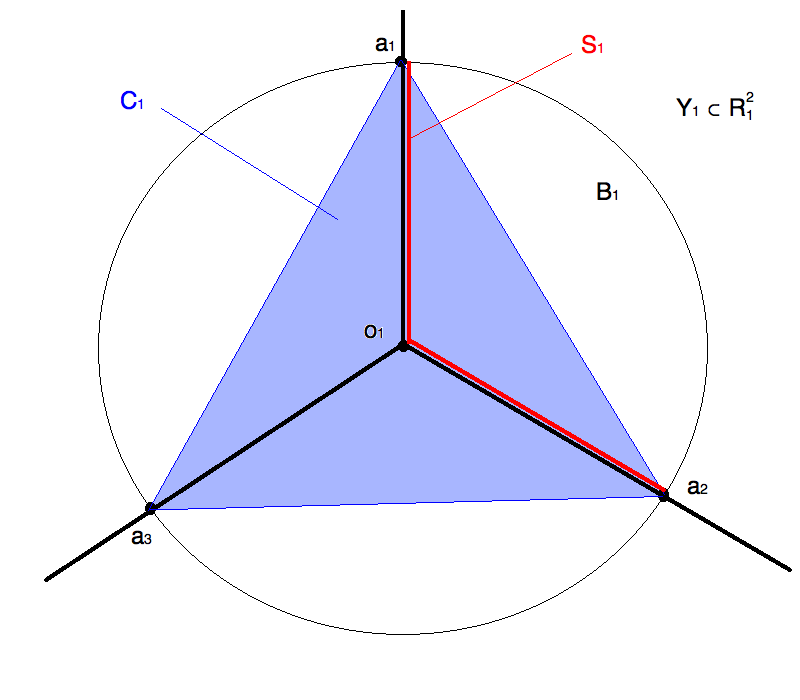}}
\nopagebreak[4]
\centerline{Figure 1} 

For $1\le j,l\le 3$, denote by $\gamma_{jl}$ the boundary of the "square" $S_j\times R_l\subset Y_1\times Y_2=Y\times Y$. Then $\gamma_{jl}\subset \partial D$, and  it represents a zero element in the singular homology group $H_1(Y\times Y\cap D, \Z/2\Z)$ of the space $Y\times Y\cap D$, because it is the boundary of the topological 2-surface $S_j\times R_l$. 

Since the set $Y\times Y$ is a cone, to prove its minimality in any sense, it suffices to consider only those competitors that only differ from $Y\times Y$ in $D$. So let us first define the fundamental class of competitors.

For any $r\in \R$, denote by $\d_r$ the dilatation map $\d_r:\R^4\to \R^4$, $\d_r(x)=rx$. For any set $E\subset \R^n$, denote by $rE=\{rx;x\in E\}$.

\begin{defn}A fundamental competitor $F$ for $Y\times Y$ is a closed set $F\subset \R^4$ such that $F\bs D^\circ=Y\times Y\bs D^\circ$ ($D^\circ$ denotes the interior of $D$), and either of the following two equivalent conditions is true:

$1^\circ$ For every $1\le j,l\le 3$, $\gamma_{jl}$ represents a zero element in the singular homology group $H_1(F\cap D, \Z/2\Z)$;

$2^\circ$ For every $1\le j,l\le 3$, and every $r\ge 1$, $r\gamma_{jl}\subset F$ represents a zero element in the singular homology group $H_1(F\cap rD, \Z/2\Z)$.

Denote by $\F$ the class of all fundamental competitors defined above.
\end{defn}

\begin{rem}Let us check the equivalence of the two conditions in Definition 1.1. 

$2^\circ\Rightarrow1^\circ$ is trivial;

$1^\circ\Rightarrow 2^\circ$: since $\gamma_{jl}$ represents a zero element in $H_1(F\cap D,\Z/2\Z)$, it represents of course a zero element in $H_1(F\cap rD,\Z/2\Z)$, since $F\cap D\subset F\cap rD$. Notice that $F$ coincides with the cone $Y\times Y$ outside $D^\circ$, hence when $r\ge 1$, for any $x\in \gamma_{jl}$, $[rt+(1-t)]x\in F\bs D^\circ$, hence the map $H_t(x)=[rt+(1-t)]x$ is a homotopy between $r\gamma_{jl}$ and $\gamma_{jl}$ on $F$. Hence the two cycles $r\gamma_{jl}$ and $\gamma_{jl}$ represent the same element in the homology group $H_1(F\cap rD,\Z/2\Z)$. As a result, $r\gamma_{jl}$ represents a zero element in $H_1(F\cap rD,\Z/2\Z)$.
\end{rem}

The singular homology group is somehow inconvenient to control, so we would like to replace it by the simplicial homology. However the simplicial homology is only defined on spaces that admit triangulations. Hence let us give the following definition of regular fundamental competitors.

\begin{defn}A set $F\in \F$ is said to be regular if there exists a finite smooth triangulation of $D$ such that $F\cap D$ is the support of a simplicial sub-complex of this triangulation. Denote by $\F^R$ the class of regular fundamental competitor. \end{defn}

The class $\F^R$ is obviously smaller than $\F$, but the following lemma will permit us to consider only regular competitors.

\begin{pro} \be \inf\{H^2(F\cap D); F\in \F\}=\inf\{H^2(F\cap D); F\in \F^R\}.\ee
\end{pro}

In the proof of Proposition 1.4, we will use mainly the construction in \cite{Fv} to construct a uniform round complex of polygons, and a Federer Fleming projection. The general purpose of the construction is to replace a set $F\in \F$ with a finite union of polygons, without increasing too much the measure. For this purpose we have to introduce first some useful notation.

In all that follows, polygons are all convex.

\begin{defn}[Polyhedral complex] Let $0<k\le n$ be integers. Let $G=\{\sigma_1,\sigma_2,\cdots,\sigma_l\}$ be a finite family of $k$-dimensional polygons of $\R^n$. For $0\le i\le k$ and $\sigma\in G$ denote by $\K_i(\sigma)$ the set of all $i$-dimensional faces of the polygon $\sigma$. Set $\K_i(G)=\cup_{\sigma\in G}\K_i(\sigma)$. Note that $\K_k(G)=G$.

Set $\K(G)=\cup_{i=0}^k\K_i(G)$.

We say that a family $\K$ of polygons of dimension at most $k$ is a polyhedral complex of dimension $k$ if there exists a family $G$ of $k-$dimensional polygons such that $\K=\K(G)$, and 
\be\forall \a,\beta\in \K, \a\ne\beta\Rightarrow \a^\circ\cap\beta^\circ=\emptyset.\ee

For $0\le i\le k$, denote by $\K_i$ the set of all its $i-$dimensional polygons, and $\K^i=\cup_{j=0}^i \K_j$ the $i-$dimensional sub-complex of $\K$.
Denote by $|\K|$ the support of the polyhedral complex $\K$:
\be |\K|=\bigcup_{\sigma\in \K}\sigma.\ee
The support $|\K^i|$ of the $i-$dimensional sub-complex $\K^i$ is called the $i$-skeleton of $\K$.
\end{defn}

\begin{defn}[Roundness constant of a polygon]Let $\sigma\in \R^n$ be a $k-$dimensional polygon, denote by $H$ the smallest affine subspace of $\R^n$ that contains $\sigma$. Then the roundness constant of $\sigma$ is 
\be R(\sigma)=1\mbox{ if }k=0;\ee
\be R(\sigma):=\frac{\inf\{r>0:\exists c\in \R^4, B(c,r)\supset\sigma\}}{\sup\{r>0:\exists c\in\sigma, B(c,r)\cap H\subset\sigma\} }\in]0,1],\mbox{ for }k\ge 1.\ee\end{defn}

\begin{rem}Roughly speaking, the roundness constant $R(\sigma)$ of a polygon $\sigma$ measures how close it is to a ball of the same dimension. The more $R(\sigma)$ is close to 1, the more it is round.\end{rem}

\begin{defn}[Roundness constant of a polyhedral complex] The roundness constant of a polyhedral complex $\K$ is
\be R(\K):=\inf_{\sigma\in \K}R(\sigma).\ee
\end{defn}

Now we want to see how to project a given closed set onto faces of polygons.

\begin{defn}[Radial projection] Let $\sigma$ be a $k-$dimensional polygon in $\R^n$, and $x\in \sigma^\circ$. Define the radial projection $\Pi_{\sigma,x}$ on the faces of $\sigma$ as follows:
\be\Pi_{\sigma, x}:=\left\{\begin{array}{rcl}\sigma\bs\{x\}&\rightarrow&\partial\sigma;\\
y&\mapsto&z\in[x,y)\cap\partial\sigma.\end{array}\right.\ee
\end{defn}

\begin{rem}Any radial projection on the faces of $\sigma$ fixes the points of $\partial \sigma$.
\end{rem}

The relation between the roundness of a polygon and the radial projections on it is given by the following lemma.
\begin{lem}[cf.\cite{Fv} Lemma 4.3.2] Let $1\le d<k\le n$ be integers. There exists a constant $K=K(d,k)>0$ that only depends on $d$ and $k$, such that for any $k-$dimensional polygon $\sigma\in \R^n$, and any set $E\subset \sigma$ with locally finite $d-$dimensional Hausdorff measure in $\sigma$, we can find a subset $X$ of $\sigma^\circ$ with non zero $H^k$ measure, such that
\be \forall x\in X, H^d(\Pi_{\sigma,x}(E))\le KR(\sigma)^{-2d}H^d(E).\ee
\end{lem}

\begin{rem}It is easy to see that $X\cap E=\emptyset$. Hence if $E$ is closed (and hence compact, because $\sigma$ is compact), then for any $x\in X$, the projection $\Pi_{\sigma,x}$ is Lipschitz on $E$ (but the Lipschitz constant could be very large).
\end{rem}

By Lemma 1.18, for $d<k\le n$, for each $k-$dimensional polyhedral complex $\K$ of roundness $R(\K)$, if $E\subset|\K|$ is a closed set with locally finite $d-$dimensional Hausdorff measure, then for each $k-$dimensional polygon $\sigma\in \K_k$, there exists a radial projection $\Pi_\sigma$ on faces of $\sigma$ such that
\be H^d(\Pi_\sigma(E\cap\sigma))\le K(d,k)R(\K)^{-2d}H^d(E\cap\sigma).\ee
Then we can define $\phi_{k-1}:E\to |\K^{k-1}|$, such that
\be \phi_{k-1}|_{\sigma}=\Pi_\sigma\mbox{ for all }\sigma\in \K_k.\ee
$\phi_{k-1}$ is well defined, because when two polygons $\a,\beta$ of the same dimension meet each other, (1.7) says that they can only meet each other at their boundaries. But  $\Pi_\a$ and $\Pi_\beta$ are both equal to the identity on boundaries, hence they agree on $\a\cap\beta$.

Set $E_{k-1}=\phi_{k-1}(E)\subset|\K^{k-1}|$. Then by (1.21) we have
\be H^d(\phi_{k-1}(E))\le K(d,k)R(\K)^{-2d}H^d(E).\ee

Now if $d=k-1$ we stop; otherwise in the $k-1$-dimensional complex $\K^{k-1}$, we can do the same thing for the $d-$dimensional subset $\phi_{k-1}(E)$ of $|\K^{k-1}|$, with a Lipschitz map $\phi_{k-2}:\phi_{k-1}(E)\to |\K^{k-2}|$ such that
\be H^d(\phi_{k-2}\circ\phi_{k-1}(E))\le K(d,k)K(d,k-1)R(\K)^{-4d}H^d(E).\ee

We carry on this process until the map $\phi_d:\phi_{d+1}\circ\cdots\circ\phi_{k-1}(E)\to |\K^d|$ is defined, with
\be H^d(\phi_d\circ\cdots\phi_{k-2}\circ\phi_{k-1}(E))\le K(d,k)K(d,k-1)\cdots K(d,d+1)R(\K)^{-2d(k-d)}H^d(E).\ee

Set $\phi'=\phi_d\circ\cdots\phi_{k-2}\circ\phi_{k-1}: E\to |\K^d|$. It is Lipschitz, and $\phi'|_{\K^d}=Id$. Set $K_1(d,k)=K(d,k)K(d,k-1)\cdots K(d,d+1)$. Then we have
\be H^d(\phi'(E))\le K_1(d,k)R(\K)^{-2d(k-d)}H^d(E).\ee

Such a $\phi'$ is called a \textbf{radial projection} (for $d-$dimensional sets) on a polyhedral complex.
  
  But we do not stop here. We want to construct a Lipschitz map $\phi:E\to|\K^d|$, such that the image $\phi(E)$ is a union of $d-$faces of $\K$, that is, if $\sigma\in \K_d$, then 
  \be \sigma^\circ\cap \phi(E)\ne\emptyset\Rightarrow \sigma\subset\phi(E).\ee
 Here for our map $\phi'$, the image $\phi'(E)$ may meet the interior of a $d$ face $\sigma$ of $\K$ without containing it. To deal with this issue, for each $\sigma\in \K_d$ that does not satisfy (1.27) with the set $\phi'(E)$, take $x\in \sigma^\circ\bs E$, and denote by $\Pi_\sigma=\Pi_{\sigma,x}$. Then $\Pi_\sigma$ is Lipschitz on $\phi'(E)\cap\sigma$ (since $E$ is compact), and it sends $\phi'(E)$ to the boundary of $\sigma$, which is of dimension $d-1$. In other words, when $\phi'(E)$ does not cover the whole $\sigma$, we "clean" it out of $\sigma$ with $\Pi_\sigma$. 
 
 Define $\phi'':\phi'(E)\to|\K^d|$: for $\sigma\in\K_d$ that satisfies (1.27) with the set $\phi'(E)$, $\phi''|_\sigma=Id$, and for $\sigma\in \K_d$ that does not satisfy (1.27) with the set $\phi'(E)$, set $\phi''|_\sigma=\Pi_\sigma$. Then $\phi'':\phi'(E)\to |\K^d|$ is 1-Lipschitz, hence 
 \be H^d(\phi''(\phi'(E))\le H^d(\phi'(E)). \ee
 
 Such a $\phi''$ is called a \textbf{polyhedral erosion}.
 
Now set $\phi=\phi''\circ\phi'$. Then $\phi$ is a Lipschitz map from $E$ to $|\K^d|$ that satisfies (1.27), and
 \be H^d(\phi(E))\le K_1(d,k)R(\K)^{-2d(k-d)}H^d(E).\ee
 
 Such a projection $\phi$ is a \textbf{Federer-Fleming projection}.
 
 Thus we have the following
 
\begin{lem}[Federer-Fleming projection]Let $1\le d<k\le n$ be integers, then there exists a constant $K_1(d,k)$ that only depends on $d$ and $k$, such that the following is true: If $\K$ is a $k-$dimensional polyhedral complex of roundness $R(\K)$, and $E\subset|\K|$ is a closed set with locally finite $d-$dimensional Hausdorff measure, then there exists Lipschitz maps $\phi',\phi''$ and $\phi$ such that

 $1^\circ$ $\phi':E\to |\K^d|$ is a radial projection, $\phi'|_{\K^d}=Id$, and satisfies (1.26);
 
 $2^\circ$ $\phi'':\phi'(E)\to |\K^d|$ is a polyhedral erosion, hence is 1-Lipschitz;
 
 $3^\circ$ $\phi=\phi''\circ\phi':E\to |\K^d|$ is a Federer Fleming projection that satisfies (1.27) and (1.29).
\end{lem}

Now we are ready to prove Proposition 1.4.

\noindent \textbf{Proof of Proposition 1.4.} We want to show that for any $F\in \F$, and any $\e>0$, there exists $F'\in \F^R$ such that 
\be H^2(F'\cap D)<H^2(F\cap D)+\e,\ee
which will give (1.5).  The main idea is to find a proper $4-$dimensional polyhedral complex whose support is $D$, and a Federer-Fleming projection $g$, such that $F'=g(F)$ verifies (1.31). The set $F'$ is clearly regular, since it is the union of $2-$faces of the polyhedral complex. In addition, we will prove that $F'$ is also a fundamental competitor.

Let $F\in\F$ be any fundamental competitor of $Y\times Y$. Without loss of generality, we can suppose that $F$ has locally finite 2-Hausdorff measure (i.e. (0.7) is true where we replace $E$ by $F$), otherwise (1.31) is automatically true. Then for any $r\in ]0,1]$, $rF\bs rD^\circ=Y\times Y\bs rD^\circ$. In particular, $rF\bs D^\circ=Y\times Y\bs D^\circ$. On the other hand, the map $\d_r: F\cap \frac1r D\to rF\cap D$ induces a homomophism between homology groups ${\d_r}_*: H_1(F\cap\frac1r D,\Z/2\Z)\to H_1(rF\cap D,\Z/2\Z)$. Now $F$ is a fundamental competitor, hence $\frac 1r\gamma_{jl}$ represents a zero element in the group $H_1(F\cap \frac1r D,\Z/2\Z)$, by Definition 1.1-$2^\circ$.  But for each $1\le j,l\le 3$, $\d_r$ sends $\frac 1r\gamma_{jl}$ to $\gamma_{jl}$, hence $\gamma_{jl}$ represents also a zero element in the group $H_1(rF\cap D,\Z/2\Z)$. This means the set $rF$ verifies $1^\circ$ in Definition 1.1.

Hence $rF$ is also a fundamental competitor for any $r\in ]0,1]$. (We will take the parameter $r$ to be very close to 1 later).

Set $\overline U_1(r):=D\bs (1-\frac{1-r}{10})D^\circ$ and $U_2(r)=(r+\frac{1-r}{10})D^\circ$. Then dist$(U_2(r),U_1(r))\ge \frac {1-r}{3}$, and $D\bs (\overline U_1(r)\cup U_2(r))$ is an annulus region. The set $rF\cap rD$ is closed in $U_2(r)$.

Now by Theorem 4.3.4 of \cite{Fv}, we can find an $n-$dimensional polyhedral complex (which is surely the support of a smooth simplicial sub-complex of a triangulation of $\R^n$) $S_r$, such that all the polyhedrons in $S_r$ are of similar sizes $\eta r$ with $\eta<\frac{1-r}{1000}$, and a deformation $\varphi$ in $U_2(r)$ (see Definition 0.1), such that

1) The support of the polyhedral complex satisfies $rD\subset |S_r|^\circ, |S_r|\subset U_2(r)$. And for each $\sigma\in S_r$, then the roundness constant $R(\sigma)$ of $\sigma$ is larger than a constant $R\in]0,1]$, that does not depend on $r$ or $rF$ or $\eta$ when $\eta$ is small. (In fact $R=R(4,2)$ depend only on the ambient dimension, which is 4, and the dimension of the set, which is 2.)

2) The map $\phi=\varphi|_{|S_r|\cap rF}$ is a radial projection. 
%
%And if we denote by $|S_r|$ the support of $S_r$, then $rD\subset |S_r|^\circ, |S_r|\subset U_2(r)$. In particular $\varphi(rF\cap U_2(r))\cap |S_r|$ is a union of $2-$faces of $S_r$.

3) $H^2(\varphi(rF)\cap U_2(r))\le H^2(rF\cap U_2(r))+(1-r)$.

\begin{rem}Notice that 3) says that $\varphi$ almost does not increase the measure of $rF$. This cannot be simply deduced from the property of a radial projection stated in Lemma 1.30, which only gives a proportional control. Here the trick is we also do something on the construction of the polyhedral complex $S_r$. Given the set $rF$, we decompose it into rectifiable part $A$ and purely unrectifiable par $B$. Then we construct the complex $S_r$ such that any of its 2-faces near $A$ are always almost parallel to the nearby tangent plane of $A$, so that the projection from $A$ to these faces almost does not increase the measure. These can be first done locally, and then be regrouped by a covering theorem. For the part $B$, the projection maps it to a set of measure 0, which is possible since it is purely unrectifiable.  

See \cite{Fv} for more detail.
\end{rem}

Next, we divide the compact region $\overline U_1(r):=D\bs (1-\frac{1-r}{10})D^\circ$ into polyhedral complex $P_r$, with roundness constant $R(P_r)$ larger than $R$, such that $Y\times Y\cap \overline U_1(r)$ is a union of 2-faces of these polygons. Moreover the sizes of the polygons in $P_r$ are about $\eta r$.

Now the support of $S_r$ is relatively far away from the cubes in $P_r$, hence by Theorem 2.3 of \cite{Fv}, there exists a polyhedral complex $S_r'$, such that $S'_r$ contains all the polygons in $S_r$ and $P_r$, $|S'_r|=D$, and the roundness constant $R(S_r')= R(P_r\cup S_r)\ge R$. (The idea is that we can carefully divide the region $D\bs (|P_r|\cup |S_r|)$ into polygons with roundness no less than $R$, and put all these polygons, and the polygons in $P_r$ and $S_r$ together to get the largest complex $S'_r$. In other words, we do some fusion to join the two complexes $P_r$ and $S_r$ together, without losing any roundness.) 

Denote by $\K_r:=(S'_r\bs S_r)\cup(S'_r\cap S_r)$ the 4-dimensional polyhedral complex. Then $|\K_r|=|S'_r|\bs |S_r|^\circ$. 

Recall that now we have four polyhedral complexes, and the relation between them is: $P_r$ and $S_r$ are far from each other. $P_r$ is near the boundary of $D$, and the support $|S_r|$ of $S_r$ contains $rD$. $S'_r\supset P_r\cup S_r$ is the biggest, whose support is $D$. $\K_r$ is almost $S'_r\bs S_r$, $|\K_r|=|S'_r|\bs|S_r|^\circ=D\bs|S_r|^\circ$.

Recall that our goal is to deform $rF$ into a regular competitor, whose measure is nearly the same as that of $rF$. In the main part $rD=|S_r|$, we have deform the set to the 2-skeleton of the polyhedral complex $S_r$, without losing much measure. Now we will deal with the remaining part of $D$, i.e. $|\K_r|$. Here when $r$ is close to $1$, the measure of $rF\cap |\K_r|$ is also small, so we just have to find a deformation that does not increase enormously its measure.

We apply Lemma 1.30 (with $k=4,d=2$) to $\K_r$, and get a  radial projection $\psi: rF\cap |\K_r|\to |\K_r^2|$ (recall that $\K_r^2$ denotes the 2-dimensional sub complex of $\K_r$), such that 
\be H^2(\psi(rF\cap|\K_r|))\le K_1(2,4)R^{-8}H^2(rF\cap|\K_r|).\ee 

Notice that near the boundary of $D$ (i.e. in $\overline U_1(r)$), $rF$ coincides with $Y\times Y$, which is a union of 2-faces of polygons of $P_r$, and hence a union of 2-faces of $\K_r$. Therefore the radial projection $\psi$ is equal to the identity on $rF\cap \overline U_1(r)$.

Now we define $f: rF\to \R^4$ as the following :
\be f(x)=\left\{\begin{array}{cc}\varphi(x), &x\in |S_r|;\\ \psi(x), &x\in |S'_r|\bs |S_r|^\circ;\\ x, &x\in \R^4\bs D.
 \end{array}\right.\ee
 Notice that $\varphi$ and $\psi$ coincides on $\partial |S_r|$ (because they are both radial projections, and hence are both identity at their boundaries), and $\psi$ coincides with $Id$ at the boundary of $D=|S'_r|$, hence the deformation $f$ is well defined, and is Lipschitz (but we do not care about the Lipschitz constant). 
  
 We are almost done. Here the image $f(rF)\cap D$ is contained in the 2-skeleton of $S'_r$. But we still want it to be regular, so we do a last step, that is, the erosion part. We apply Lemma (1.30) to $f(rF)\cap D$, and get a polyhedral erosion $h:f(rF)\cap D\to|\K_r^2|$, which is 1-Lipschitz, and the image $h(f(rF)\cap D)$ is a union of 2-faces of $\K_r$, and hence is regular.
 
 Still notice that at the boundary of $D$, $f(rF)$ is already a union of $2$-faces of $\K_r$, hence $h|_{\partial D}=Id$. So we can safely define the Lipschitz map
 \be \begin{split}g&:rF\to \R^4; \\g(x)&=\left\{\begin{array}{rcl}h\circ f(x)&, &x\in rF\cap D; \\f(x)&, &x\in rF\bs D.\end{array}\right.\end{split}\ee
 
The map $g$ is Lipschitz, thus we can extend $g$ to a Lipschitz map from $\R^4$ to $\R^4$, with $g|_{\R^4\bs D}=Id$. We still call this extension $g$ for short. Notice also that $\varphi,\psi,f$ and $g$ depend on $r$. But the constants $R$ and $K_1(2,4)$ in (1.33) do not.

Now let us look at the set $F_r=g(rF)$. Notice first that in $\overline U_1(r):=D\bs (1-\frac{1-r}{10})D^\circ$, $rF$ coincides with $Y\times Y$, which is a union of 2-faces of polygons of $P_r$ and hence of $S'_r$. Hence $g(rF)\bs (1-\frac{1-r}{10})D^\circ=Y\times Y\bs (1-\frac{1-r}{10})D^\circ$. This gives that $F_r\cap\partial D=Y\times Y\cap\partial D$. In particular, $F_r\bs D^\circ=Y\times Y\bs D^\circ$. On the other hand, $F_r$ is the image of the map $g$ of $rF$, with $g(rF\cap D)=F_r\cap D$, and $g$ fixes every $\gamma_{jl}, 1\le j,l\le 3$, hence the map $g$ induces a group homomorphism $g_*: H_1(rF\cap D, \Z/2\Z)\to H_1(F_r\cap D, \Z/2\Z)$. As a result, since each $\gamma_{jl}, 1\le j,l\le 3$ represents a zero element in $H_1(rF\cap D, \Z/2\Z)$ (recall that we have proved at the beginning that $rF$ is a fundamental competitor), it represents a zero element in $H_1(F_r\cap D, \Z/2\Z)$. Hence $F_r=f(rF)\in\F$. Moreover, since $F_r\cap D$ is a sub-complex of $S'_r$, $F_r$ is regular, and hence $F_r\in \F^R$.

Now the thing left is to estimate the measure of $F_r\cap D$. We have
\be \begin{split}H^2(F_r\cap D)&\le H^2(F_r\cap |S_r|)+H^2(F_r\cap |\K_r|)\\
&=H^2(g(rF)\cap |S_r|)+H^2(g(rF)\cap |\K_r|)\\
&=H^2(h\circ f(rF)\cap |S_r|)+H^2(h\circ f(rF)\cap |\K_r|)\\
&\le H^2(f(rF)\cap |S_r|)+H^2(f(rF)\cap |\K_r|) \hfill\mbox{(because }h\mbox{ is 1-Lipschitz)}\\
&=H^2(\varphi(rF)\cap U_2(r))+H^2(\psi(rF\cap |\K_r|)\\
&<H^2(rF\cap U_2(r))+(1-r)+K_1(2,4)R^{-8}H^2(rF\cap |\K_r|)\\
&<H^2(rF\cap U_2(r))+(1-r)+K_1(2,4)R^{-8}H^2(rF\cap D\bs rD)\\
&=H^2(rF\cap U_2(r))+(1-r)+K_1(2,4)R^{-8}H^2(Y\times Y\cap D\bs rD),\end{split}\ee
The last equality is because $rF\cap D\bs rD=D\cap r(F\bs D)=D\cap r(Y\times Y\bs D)=Y\times Y\cap D\bs rD$.

Notice that the two terms $(1-r)$ and $K_1(2,4)R^{-8}H^2(Y\times Y\cap D\bs rD)$ tend to zero as $r$ tends to 1, hence for each $\e>0$, there exists a $r=r_\e\in ]0,1[$ such that $H^2(F_r\cap D)<H^2(F\cap D)+\e$. Set $F'=F_r\in \F^R$, then we get (1.5). \qed

%Denote by $\F_D$ the class $\{F\in\F: F\bs D=E\bs D\}$. 

\section{Decomposition of a regular fundamental competitor in D}

In this section we will decompose regular fundamental competitors of $Y\times Y$ in $D$ into pieces that do not make inconvenient intersections.

Let $F\in\F^R$ be a regular fundamental competitor for $Y\times Y$ in $D$. It is the support of a simplicial complex, hence the singular homology groups and simplicial homology groups are naturally isomorphic (cf. \cite{Mun}, Theorem 34.3).

For $1\le j,l\le 2$, since $\gamma_{jl}$ represents a zero element in $H_1(F\cap D,\Z/2\Z)$, there exists a simplicial 2-chain $\Gamma_{jl}\subset F\cap D$ with coefficient in $\Z/2\Z$ such that $\partial \Gamma_{jl}=\gamma_{jl}.$

Now denote by $\Gamma_{13}$ the simplicial 2-chain $\Gamma_{13}=\Gamma_{11}+\Gamma_{12}$ with coefficients in $\Z/2\Z$. Let us check that $\partial\Gamma_{13}=\partial\Gamma_{11}+\partial\Gamma_{12}=\gamma_{11}+\gamma_{12}=\gamma_{13}$. In fact, by definition of $\gamma_{ij}$, (notice that the coefficient group is $\Z/2\Z$)
\be \begin{split}\gamma_{11}&=\{a_1\}\times R_1+\{a_2\}\times R_1+\{b_1\}\times S_1+\{b_2\}\times S_1,\\
\gamma_{12}&=\{a_1\}\times R_2+\{a_2\}\times R_2+\{b_2\}\times S_1+\{b_3\}\times S_1,\end{split}\ee
hence
\be \begin{split}\gamma_{11}+\gamma_{12}&=\{a_1\}\times (R_1+R_2)+\{a_2\}\times (R_1+R_2)+\{b_1\}\times S_1+\{b_3\}\times S_1\\
&=\{a_1\}\times R_3+\{a_2\}\times R_3+\{b_1\}\times S_1+\{b_3\}\times S_1=\gamma_{13}.\end{split}\ee

It is important to point out that, by definition of $\Gamma_{13}$, modulo a $H^2$ negligible set, the supports of the three 2-chains satisfy that :
\be |\Gamma_{13}|=(|\Gamma_{12}|\cup|\Gamma_{11}|)\bs(|\Gamma_{12}|\cap|\Gamma_{11}|).\ee
In other words, $H^2-$almost every point in the union $|\Gamma_{12}|\cup|\Gamma_{11}|\cup|\Gamma_{13}|$ of the three supports  belongs to exactly two of them.

We define in the same way :
\be \Gamma_{23}=\Gamma_{21}+\Gamma_{22}.\ee
Then by the same argument we have
\be \partial\Gamma_{23}=\gamma_{23}\mbox{ and }|\Gamma_{23}|=(|\Gamma_{22}|\cup|\Gamma_{21}|)\bs(|\Gamma_{22}|\cap|\Gamma_{21}|)\ee
modulo a $H^2-$negligible set.

Thus we have now defined the $\Gamma_{jl}$ for $1\le j\le 2$ and $1\le l\le3.$ We continue to define, for $1\le l\le 3$, 
\be \Gamma_{3l}=\Gamma_{1l}+\Gamma_{2l}.\ee
Then by the same principle, we have
\be \partial\Gamma_{3l}=\gamma_{1l}\mbox{ and }|\Gamma_{2l}|=(|\Gamma_{1l}|\cup|\Gamma_{2l}|)\bs(|\Gamma_{1l}|\cap|\Gamma_{2l}|)\ee
modulo a $H^2-$negligible set.

The last thing to check is that modulo a $H^2-$negligible set, 
\be |\Gamma_{33}|=(|\Gamma_{31}|\cup|\Gamma_{32}|)\bs(|\Gamma_{31}|\cap|\Gamma_{32}|).\ee

But this follows directly from the definition, since
\be\begin{split} 
\Gamma_{33}&=\Gamma_{13}+\Gamma_{23}=\Gamma_{11}+\Gamma_{12}+\Gamma_{21}+\Gamma_{22}\\
&=(\Gamma_{11}+\Gamma_{21})+(\Gamma_{12}+\Gamma_{22})=\Gamma_{31}+\Gamma_{32}.
\end{split}\ee

 To sum up, we have these nine subsets $|\Gamma_{jl}|, 1\le j,l\le 3$, of our regular fundamental competitor $F$ of $Y\times Y$, which are supports of smooth simplicial 2-chains in $\R^4$, such that for any $1\le j\le 3$,  $H^2-$ almost every point in the union of the three supports $|\Gamma_{j1}|,|\Gamma_{j2}|,|\Gamma_{j3}|$ belongs 
to exactly two of them. The same holds for $|\Gamma_{1j}|,|\Gamma_{2j}|,|\Gamma_{3j}|$.

\bigskip

Now we have to prove some projection property for these nine subsets of $F$. Denote by $L_1$ the segment connecting $a_1$ and $a_2$, $L_2$ the segment connecting $a_2$ and $a_3$, $L_3$ the segment connecting $a_3$ and $a_1$; $M_1$ the segment connecting $b_1$ and $b_2$, $M_2$ the segment connecting $b_2$ and $b_3$, and $M_3$ the segment connecting $b_3$ and $b_1$. Denote by $Q_{jl}, 1\le j,l\le 3$ the square $L_j\times M_l$. Denote by $P_{jl},1\le j,l\le 3$ the plane containing $Q_{jl}$. Denote by $p_{jl}$ the orthogonal projection from $\R^4$ to $P_{jl}$.

\begin{lem}For $1\le j,l\le 3$, the set $|\Gamma_{jl}|$ satisfies
\be p_{jl}(|\Gamma_{jl}|)\supset Q_{jl}.\ee
\end{lem}
 
 \nd Notice that $p_{jl}: |\Gamma_{jl}|\to P_{jl}$ induces a group homomorphism ${p_{jl}}_*: H_1(|\Gamma_{jl}|,\Z/2\Z)\to H_1(p_{jl}(|\Gamma_{jl}|),\Z/2\Z)$, where $p_{jl}(|\Gamma_{jl}|)\subset Q_{jl}$. Hence the image of the chain ${p_{jl}}_*(\gamma_{ij})=\partial Q_{jl}$ is a zero element in $H_1(p_{jl}(|\Gamma_{jl}|),\Z/2\Z)$. Now if the projection $p_{jl}(|\Gamma_{jl}|)$ does not contain $Q_{jl}$, for example there exists $x\in Q_{jl}\bs p_{jl}(|\Gamma_{jl}|)$, then a radial projection from $Q_{jl}\bs\{x\}\to\partial Q_{jl}$ will map homotopically $p_{jl}(|\Gamma_{jl}|)$ to $\partial Q_{jl}$. Hence we have $\partial Q_{jl}$ represents also a zero element in $H_1(\partial Q_{jl},\Z/2\Z)$, this gives  a contradiction.\qed
 
 Now for convenient use in the next section, denote $F_{jl}=|\Gamma_{jl}|, 1\le j,l\le 3$. Then to sum up, we proved the following proposition of this section:
 
 \begin{pro}For every regular fundamental competitor $F$ of $Y\times Y$ in $D$, there exist subsets $F_{jl},1\le j,l\le 3$ of $F$, which are supports of smooth simplicial 2-chains $\Gamma_{jl}$ in $\R^4$, such that
 \be p_{jl}(F_{jl})\supset Q_{jl},\ee
 and for any $1\le j\le 3$,  $H^2-$almost every point in the union of the three subsets $F_{j1},F_{j2},F_{j3}$ (resp. $F_{1j},F_{2j},F_{3j}$) belongs 
to exactly two of them. 
 \end{pro}
 
 \section{The calibration}
 
In this section we will use a product of paired calibrations to prove the first theorem of this paper:

\begin{thm}\be H^2(Y\times Y\cap D)=\inf\{H^2(F\cap D); F\in \F\}.\ee
\end{thm}

\nd By Proposition 1.4, we just have to prove that
\be H^2(Y\times Y\cap D)=\inf\{H^2(F\cap D); F\in \F^R\}.\ee

We first define the calibration.

For each $1\le j\le 3$, denote by $x_j$ the  midpoint of the segment $L_j$, and set $w_j:=\frac{\overrightarrow{o_1x_j}}{|\overrightarrow{o_1x_j}|}$, which is a unit vector in $\R^2_1$ that is normal to $L_j$. Similarly, for $1\le l\le 3$ we denote by $y_l$ the  midpoint of the segment $M_l$, and a unit vector $u_l:=\frac{\overrightarrow{o_2y_l}}{|\overrightarrow{o_2y_l}|}$ in $\R^2_2$ that is normal to $M_l$. Then, for $1\le j,l\le 3$, denote by $v_{jl}=w_j\wedge u_l$, a unit simple 2-vector in $\R^4$.

We define the function $f_{jl}$ on the set of simple 2-vectors in $\R^4$ : for any simple 2-vector $\xi\in \wedge_2(\R^4)$, $f_{jl}(\xi):=|\xi\wedge v_{jl}|=|\det_{e_1,e_2,e_3,e_4}\xi\wedge v_{jl}|$, with $\{e_j\}_{1\le j\le 4}$ the canonical orthonormal basis of $\R^4$. Now for any unit (with respect to the $L^2$ norm $|\cdot|$ for the orthonormal basis $\{e_i\wedge e_j\}_{1\le i<j\le 4}$ of $\wedge_2(\R^4)$) simple 2-vector $\xi$, we can associate to it a plane $P(\xi)\in G(4,2)$, where $G(4,2)$ is the set of all 2-dimensional subspaces of $\R^4$ :
\be P(\xi)=\{v\in \R^4, v\wedge\xi=0\}.\ee
In other words, $P(x\wedge y)$ is the subspace generated by $x$ and $y$.

Now denote also by $g_{jl}$ the function from $G(4,2)$ to $\R$: for any $P=P(x\wedge y)\in G(4,2)$ with $x\wedge y$ a unit simple 2-vector, $g_{jl}(P)=f_{jl}(x\wedge y)$. Since the definition of $f_{jl}$ on $\wedge_2(\R^4)$ is to take the absolute value of the determinant, the function $g_{jl}$ is well defined.

Let $F_{jl}, 1\le j,l\le 3$ be as in Proposition 2.12. Now since $F_{jl}$ is the support of a smooth simplicial 2-chain in $\R^4$, it is 2-rectifiable., and the tangent plane $T_xF_{jl}$ of $F_{jl}$ at  $x$ exists for $H^2-$almost all $x\in F_{jl}$. We want to estimate 
\be \int_{F_{jl}}g_{jl}(T_xF_{ij}) dH^2(x).\ee

Notice that $F_{ij}$ is piecewise smooth, hence $g_{jl}(T_xF_{ij})$ is measurable. Note also that $||g_{jl}||_\infty=1$, hence the integral is well defined.

Denote by $E_{jl}:=\{x\in F_{jl} : J_2p_{jl}(x)\ne 0\}\subset F_{ij}$, where $J_2p_{jl}(x)=||\wedge_2D(p_{jl}|_{F_{jl}})(x)||$ is the Jacobian of the restriction $p_{jl}|_{F_{jl}}$. Then $T_xE_{jl}=T_xF_{jl}$ for $H^2$ almost all $x\in E_{jl}$. By the Sard theorem, we have
\be H^2(p_{jl}(F_{jl})\bs p_{jl}(E_{jl}))=0,\ee
and hence by Proposition 2.12, 
\be H^2(p_{jl}(E_{jl}))=H^2(p_{jl}(F_{jl}))\ge H^2(Q_{jl}).\ee
 
 Now for $x\in E_{jl}$, define $h_{jl}(x)=g_{jl}(T_xF_{jl}) (J_2p_{jl}(x))^{-1}$. 
Recall that the projection $p_{jl}: F_{jl}\to P_{jl}$ is a 1-Lipschitz function.  hence by the coarea formula for Lipschitz functions between two rectifiable set (cf.\cite{Fe} Theorem 3.2.22), we have
\be \int_{E_{jl}}h_{jl}(x) J_2p_{jl}(x) dH^2(x)=\int_{p_{jl}(E_{jl})} dH^2(y)[\sum_{p_{jl}(x)=y}h_{jl}(x)].\ee

By definition of $h_{jl}$, the left hand side of the above equality is just
\be \begin{split}\int_{E_{jl}}h_{jl}(x) J_2p_{jl}(x) dH^2(x)&=\int_{E_{jl}}g_{jl}(T_xF_{jl}) (J_2p_{jl}(x))^{-1} J_2p_{jl}(x) dH^2(x)\\
&=\int_{E_{jl}}g_{jl}(T_xF_{jl})dH^2(x)\le \int_{F_{jl}}g_{jl}(T_xF_{jl})dH^2(x),\end{split}\ee
the last inequality is because $g_{jl}$ is non negative.

For the right hand side, note that for any $x\in E_{jl}$, suppose $T_xF_{jl}=P(u\wedge v)$, with $u,v$ an orthonormal basis of $T_x(F_{ij})$. Hence
\be g_{jl}(T_xF_{jl})=f_{jl}(u\wedge v)=|v_{jl}\wedge u\wedge v|=|v_{jl}\wedge p_{jl}(u\wedge v)|.\ee

Notice that $p_{jl}(u\wedge v)\in\wedge_2(P_{jl})$, hence if we take a unit simple two vector $\xi_{jl}$ of $P_{jl}$, we have $p_{jl}(u\wedge v)=\pm|p_{jl}(u\wedge v)|\xi_{jl},$ and hence by (3.10)
\be g_{jl}(T_xF_{jl})=|p_{jl}(u\wedge v)||v_{jl}\wedge\xi_{jl}|=|p_{jl}(u\wedge v)|=J_2p_{jl}(x),\ee
and thus 
\be h_{jl}(x)=g_{jl}(T_xF_{jl})J_2p_{jl}(x)^{-1}=1\ee
 for $H^2-a.e.$ $x\in E_{jl}$.
 As a result, 
 \be \int_{p_{jl}(E_{jl})} dH^2(y)[\sum_{p_{jl}(x)=y}h_{jl}(x)]=\int_{p_{jl}(E_{jl})} dH^2(y)[\sum_{p_{jl}(x)=y}1]
 =\int_{p_{jl}(E_{jl})} dH^2(y)\sharp \{p_{jl}^{-1}(x)\cap E_{jl}\}.\ee
But for all $y\in p_{jl}(E_{jl})$, $\sharp \{p_{jl}^{-1}(x)\cap E_{jl}\}\ge 1$, hence
\be \int_{p_{jl}(E_{jl})} dH^2(y)[\sum_{p_{jl}(x)=y}h_{jl}(x)]\ge \int_{p_{jl}(E_{jl})} dH^2(y)=H^2(p_{jl}(E_{jl}))\ge H^2(Q_{jl}).\ee

Combine (3.13) (3.9) and (3.8) we have
\be  H^2(Q_{jl})\le \int_{F_{jl}}g_{jl}(T_xF_{jl})dH^2(x), \mbox{ for }1\le j,l\le 3.\ee

We sum over $1\le j,l\le 3$, and have
\be \begin{split}\sum_{1\le j\le 3}\sum_{1\le l\le 3}H^2(Q_{jl})&\le \sum_{1\le j\le 3}\sum_{1\le l\le 3}\int_{F_{jl}}g_{jl}(T_xF_{jl})dH^2(x)\\
&=\int_{\cup_{1\le j,l\le 3}F_{jl}}dH^2(x) [\sum_{1\le j\le 3}\sum_{1\le l\le 3}g_{jl}(T_xF_{jl})1_{F_{jl}}(x)].\end{split}\ee

But $F$ is 2-rectifiable, and each $F_{jl}$ is its subset, hence we have for $H^2-a.e. x\in F_{jl}$, $T_xF_{jl}=T_xF$. Hence we have
\be \sum_{1\le j\le 3}\sum_{1\le l\le 3}H^2(Q_{jl})\le \int_{\cup_{1\le j,l\le 3}F_{jl}}dH^2(x) [\sum_{1\le j\le 3}\sum_{1\le l\le 3}g_{jl}(T_xF)1_{F_{jl}}(x)].\ee

Now we want to use Proposition 2.12 to derive a essential upper bound for the function 
\be [\sum_{1\le j\le 3}\sum_{1\le l\le 3}g_{jl}(T_xF)1_{F_{jl}}(x)].\ee

Given a point $x\in \cup_{1\le j,l\le 3}F_{jl}$, then by Proposition 2.12, modulo a negligible set, there are two possibilities :

1) There exists $1\le j_1,j_2\le 3$ and $1\le l_1,l_2\le 3$ such that $x$ only belongs to the four pieces $F_{j_1l_1},F_{j_1l_2},F_{j_2l_1},F_{j_2l_2}.$

2) There exists a permutation $\sigma$ of $\{1,2,3\}$ such that $x$ belongs to all the nine $F_{jl}$ except for $F_{1\sigma(1)},F_{2\sigma(2)},F_{3\sigma(3)}$.

We will estimate the function $[\sum_{1\le j\le 3}\sum_{1\le l\le 3}g_{jl}(T_xF)1_{F_{jl}}(x)]$ in these two cases.

For 1), without loss of generality, we suppose that $j_1=l_1=1,j_2=l_2=2$.  Then 
\be \sum_{1\le j\le 3}\sum_{1\le l\le 3}g_{jl}(T_xF)1_{F_{jl}}(x)=\sum_{1\le j\le 2}\sum_{1\le l\le 2}g_{jl}(T_xF).\ee
 Suppose that $T_xF=P(\xi)$ with $\xi$ a unit simple 2-vector. Then for each $1\le j,l\le 2$, by definition of $g_{jl}$, $g_{jl}(T_xF)=|v_{jl}\wedge\xi|$.

Hence 
\be\begin{split}\sum_{1\le j\le 2}\sum_{1\le l\le 2}g_{jl}(T_xF)&=\sup_{\e}\sum_{1\le j,l\le 2}\e(j,l)\det(v_{jl}\wedge\xi)\\
&=\sup_{\e}\det[(\sum_{1\le j,l\le 2}\e(j,l)v_{jl})\wedge\xi]\le\sup_{\e}||\sum_{1\le j,l\le 2}\e(j,l)v_{jl}||.\end{split}\ee
where $\e$ run over all function from $\{1,2\}\times\{1,2\}\to\{1,-1\}$, and the norm $||\cdot ||$ on $\wedge_2(\R^4)$ is defined by 
\be ||\a||=\sup\{\det(\a\wedge \beta); \beta \mbox{ simple unit 2-vector}\}.\ee
Then the last equality is because $|\xi|=1$. Hence
\be \sum_{1\le j\le 3}\sum_{1\le l\le 3}g_{jl}(T_xF)1_{F_{jl}}(x)\le \sup\{||\sum_{1\le j,l\le 2}\e(j,l)v_{jl}||,\e:\{1,2\}\times\{1,2\}\to\{1,-1\}\}.\ee

Similarly, for 2), we have
\be \begin{split}\sum_{1\le j\le 3}&\sum_{1\le l\le 3}g_{jl}(T_xF)1_{F_{jl}}(x)=\sum_{1\le j,l\le 3, l\ne\sigma(j)}g_{jl}(T_xF)\\
&\le\sup\{||\sum_{1\le j,l\le 3, l\ne  j}\e(j,l)v_{jl}||,\e:[\{1,2,3\}\times\{1,2,3\}]\bs\{(1,1),(2,2),(3,3)\}\to\{1,-1\}\}.\end{split}\ee

The following lemma will lead to the conclusion of Theorem 3.1.

\begin{lem}\be \sup\{||\sum_{1\le j,l\le 2}\e(j,l)v_{jl}||,\e:\{1,2\}\times\{1,2\}\to\{1,-1\}\}\le 3,\ee
and
\be \sup\{||\sum_{1\le j,l\le 3, l\ne j}\e(j,l)v_{jl}||,\e:[\{1,2,3\}\times\{1,2,3\}]\bs\{(1,1),(2,2),(3,3)\}\to\{1,-1\}\}\le 3.\ee
\end{lem}

We will leave the proof of this lemma to the next section. Now let us admit this lemma and finish the proof of Theorem 3.1.

By Lemma 3.24, (3.22) and (3.23), for $H^2-$almost all $x\in \cup_{1\le j,l\le 3}F_{jl}$, 
\be\sum_{1\le j\le 3}\sum_{1\le l\le 3}g_{jl}(T_xF)1_{F_{jl}}(x)\le 3.\ee
Hence by (3.17), we have
\be \begin{split}\sum_{1\le j\le 3}\sum_{1\le l\le 3}H^2(Q_{jl})&\le \int_{\cup_{1\le j,l\le 3}F_{jl}}dH^2(x) [\sum_{1\le j\le 3}\sum_{1\le l\le 3}g_{jl}(T_xF)1_{F_{jl}}(x)]\\
&\le 3\int_{\cup_{1\le j,l\le 3}F_{jl}}dH^2(x)=3H^2(\cup_{1\le j,l\le 3}F_{jl})\le 3H^2(F\cap D) .\end{split}\ee

However, notice that $Q_{jl}=L_j\times M_l$, and the length of the segments $L_j,M_l$ are $\sqrt 3$ (cause they are the edges of regular triangles inscribed to the unit ball), hence $H^2(Q_{jl})=3$. As a result
\be \sum_{1\le j\le 3}\sum_{1\le l\le 3}H^2(Q_{jl})=27=3H^2((Y_1\cap B_1)\times(Y_2\cap B_2))=3H^2(Y\times Y\cap D).\ee
Therefore by (3.28),
\be H^2(Y\times Y\cap D)\le H^2(F\cap D).\ee

Recall that $F\in \F^R$ is an arbitrary regular competitor for $Y\times Y$ in $D$. Hence we get (3.3), and the proof of Theorem 3.1 is completed.\qed

\section{Proof of Lemma 3.24}

In this section we give the prove of Lemma 3.24. We will need the following lemma for the proof.

\begin{lem}Suppose that $\{x_i\}_{1\le i\le 4}$ is an orthonormal basis of $\R^4$, then
\be ||x_1\wedge x_2\pm x_3\wedge x_4||=1.\ee
\end{lem}

\nd Let $\xi$ be a unit simple 2-vector in $\R^4$. Then there exists two unit vectors $y,z\in \R^4$ with $y\perp z$, such that $\xi=y\wedge z$. Since $\{x_i\}_{1\le i\le 4}$ is an orthonormal basis, there exists $A,B,C,D,a,b,c,d\in \R$ such that 
\be y=Ax_1+Bx_2+Cx_3+Dx_4, z=ax_1+bx_2+cx_3+dx_4,\ee
\be A^2+B^2+C^2+D^2=a^2+b^2+c^2+d^2=1,\ee
and
\be Aa+Bb+Cc+Dd=0.\ee
Thus 
\be \begin{split}\xi=y\wedge z&=(Ab-aB)x_1\wedge x_2+(Ac-aC)x_1\wedge x_3+(Ad-aD)x_1\wedge x_4\\
&+(Bc-bC)x_2\wedge x_3+(Bd-bD)x_2\wedge x_4+(Cd-cD)x_3\wedge x_4,\end{split}\ee
and hence
\be (x_1\wedge x_2\pm x_3\wedge x_4)\wedge\xi=[\pm(Ab-aB)+ (Cd-cD)]x_1\wedge x_2\wedge x_3\wedge x_4,\ee
which yields
\be \det ((x_1\wedge x_2+x_3\wedge x_4)\wedge\xi)= [\pm(Ab-aB)+ (Cd-cD)].\ee
But by (4.4)
\be [\pm(Ab-aB)+(Cd-cD)]\le\frac12(A^2+b^2+a^2+B^2+C^2+d^2+c^2+D^2)=1,\ee
hence
\be \det ((x_1\wedge x_2\pm x_3\wedge x_4)\wedge\xi)\le 1.\ee
Therefore
\be ||x_1\wedge x_2\pm x_3\wedge x_4||\le 1\ee
by definition of $||\cdot||$.

To prove the equality, it suffices to notice that
\be (x_1\wedge x_2\pm x_3\wedge x_4)\wedge (x_3\wedge x_4)=1,\ee
where $x_3\wedge x_4$ is a unit simple 2-vector. \qed

\noindent \textbf{Proof of Lemma 3.24.}

Let us first give some notation. 

We say that two indices in $\{1,2,3\}\times \{1,2,3\}$ are adjacent if they have a common left index or right index. That is, $(j,k)$ and $(j,l)$ are adjacent, and $(k,l)$ and $(j,l)$ are adjacent. 

\textbf{Proof of (3.25). }

Case 1. The map $\e$ is such that there exists two adjacent indices $a,b\in \{1,2\}\times\{1,2\}$ with $\e(a)=\e(b)$. Without loss of generality, suppose that $\e(1,1)=\e(1,2)$. Notice that $v_{11}+v_{12}+v_{13}=0$, hence
\be \sum_{1\le j,l\le 2}\e(j,l)v_{jl}=\e_{1,1}v_{11}+\e_{1,2}v_{12}+\e_{2,1}v_{21}+\e_{2,2}v_{22}=-\e_{1,1}v_{13}+\e_{2,1}v_{21}+\e_{2,2}v_{22},\ee
which gives
\be ||\sum_{1\le j,l\le 2}\e(j,l)v_{jl}||\le ||v_{13}||+||v_{21}||+||v_{22}||.\ee

By definition of $||\cdot||$ in (3.20), for any $\a\in\wedge_2(\R^4)$, for any $\beta$ unit simple 2-vector in $\R^4$, 
\be |\det(\a\wedge\beta)|=|\a\wedge\beta|\le|\a||\beta|=|\a|,\ee
because $|\cdot|$ is the $L^2$ norm. As a result,
\be ||\a||=\sup\{\det(\a\wedge\beta);\beta\mbox{ unit simple 2-vector}\}\le |\a|.\ee
Combine this with (4.14), we get
\be ||\sum_{1\le j,l\le 2}\e(j,l)v_{jl}||\le |v_{13}|+|v_{21}|+|v_{22}|=3.\ee

Case 2. For any two adjacent indices $a,b\in\{1,2\}\times\{1,2\}$, $\e(a)=-\e(b)$. Then
\be \e(1,1)=\e(2,2)=-\e(1,2)=-\e(2,1).\ee
In this case,
\be \begin{split}\sum_{1\le j,l\le 2}\e(j,l)v_{jl}&=\e(1,1)(v_{11}-v_{12}+v_{22}-v_{21})\\
&=\e(1,1)(w_1\wedge u_1-w_1\wedge u_2-w_2\wedge u_1+w_2\wedge u_2)=\e(1,1)(w_1-w_2)\wedge (u_1-u_2),\end{split}\ee
hence
\be ||\sum_{1\le j,l\le 2}\e(j,l)v_{jl}||=||(w_1-w_2)\wedge (u_1-u_2)||\le |(w_1-w_2)\wedge (u_1-u_2)|=|(w_1-w_2)||(u_1-u_2)|=3.\ee

Combine Case 1 and Case 2, we get (3.25).

\textbf{Proof of (3.26).}

Case 1.  If $\e:[\{1,2,3\}\times\{1,2,3\}\bs\{(1,1),(2,2)(3,3)\}\to\{1,-1\}$ attributes to every pair of adjacent indices different values, then we have
\be \e(1,2)=-\e(1,3)=\e(2,3)=-\e(2,1)=\e(3,1)=-\e(3,2).\ee
Take an orthonormal basis $\{f_i\}_{1\le i\le 4}$ with $f_1,f_2\in\R^2_1,f_3,f_4\in\R^2_2$, such that
\be \begin{split}w_1=f_1,w_2=-\frac12 f_1+\frac{\sqrt 3}{2}f_2, w_3=-\frac12 f_1-\frac{\sqrt 3}{2}f_2,\\
u_2=f_3,u_3=-\frac12 f_3+\frac{\sqrt 3}{2}f_4, u_1=-\frac12 f_3-\frac{\sqrt 3}{2}f_4,\end{split}\ee
we get
\be \begin{split}&||\sum_{1\le j,l\le 3, l\ne j}\e(j,l)v_{jl}||=||v_{12}-v_{13}+v_{23}-v_{21}+v_{31}-v_{32}||\\
&=||f_1\wedge f_3-f_1\wedge (-\frac12 f_3+\frac{\sqrt 3}{2}f_4)+(-\frac12 f_1+\frac{\sqrt 3}{2}f_2)\wedge(-\frac12 f_3+\frac{\sqrt 3}{2}f_4)\\
&-(-\frac12 f_1+\frac{\sqrt 3}{2}f_2)\wedge (-\frac12 f_3-\frac{\sqrt 3}{2}f_4)+(-\frac12 f_1-\frac{\sqrt 3}{2}f_2)\wedge(-\frac12 f_3-\frac{\sqrt 3}{2}f_4)-(-\frac12 f_1-\frac{\sqrt 3}{2}f_2)\wedge f_3||\\
&=||\frac94 f_1\wedge f_3-\frac{3\sqrt 3}{4}f_1\wedge f_4+\frac{3\sqrt 3}{4}f_2\wedge f_3+\frac 94 f_2\wedge f_4||\\
&=\frac{3\sqrt 3}{2}||f_1\wedge (\frac{\sqrt 3}{2}f_3-\frac12 f_4)+f_2\wedge (\frac 12 f_3+\frac{\sqrt 3}{2} f_4)||\end{split}\ee
Notice that the two vectors $\frac{\sqrt 3}{2}f_3-\frac12 f_4,\frac 12 f_3+\frac{\sqrt 3}{2} f_4$ form an orthonormal basis of $R_2^2$, hence $\{f_1,f_2,\frac{\sqrt 3}{2}f_3-\frac12 f_4,\frac 12 f_3+\frac{\sqrt 3}{2} f_4\}$ form an orthonormal basis of $\R^4$. By lemma 4.1, $||f_1\wedge (\frac{\sqrt 3}{2}f_3-\frac12 f_4)+f_2\wedge (\frac 12 f_3+\frac{\sqrt 3}{2} f_4)||\le 1$, hence 
\be ||\sum_{1\le j,l\le 3, l\ne j}\e(j,l)v_{jl}||\le \frac{3\sqrt 3}{2}<3.\ee

Case 2. If $\e$ attribute to a pair of adjacent indices the same value. Without loss of generality, suppose $\e(1,2)=\e(1,3)$. Then
\be \e(1,2)v_{12}+\e(1,3)v_{13}=-\e(1,2)v_{11},\ee
and hence
\be\sum_{1\le j,l\le 3, l\ne j}\e(j,l)v_{jl}=-\e(1,2)v_{11}+\e(2,1)v_{21}+\e(2,3)v_{23}+\e(3,1)v_{31}+\e(3,2)v_{32}.\ee

Subcase 2.1. $\e(2,1)=\e(3,1)$, then $\e(2,1)v_{21}+\e(3,1)v_{31}=-\e(2,1)v_{11}$. There are two possibilities:

Subsubcase 2.1.1: $-\e(1,2)=\e(2,1)=\e(3,1)$. In this case, $-\e(1,2)v_{11}+\e(2,1)v_{21}+\e(3,1)v_{31}=0$, and hence
\be ||\sum_{1\le j,l\le 3, l\ne j}\e(j,l)v_{jl}||=||\e(2,3)v_{23}+\e(3,2)v_{32}||\le ||v_{23}||+||v_{32}||<3.\ee

Subsubcase 2.1.2: $\e(1,2)=\e(2,1)=\e(3,1)$. In this case we have
\be \sum_{1\le j,l\le 3, l\ne j}\e(j,l)v_{jl}=-2\e(1,2)v_{11}+\e(2,3)v_{23}+\e(3,2)v_{32}.\ee

Modulo symmetry of $v_{23}$ and $v_{32}$, and multiplication of $\pm 1$ (which does not change the norm anyway), we have three possibilities :

$1^\circ$ $2 v_{11}+v_{23}+v_{32}$;

$2^\circ$ $2v_{11}-v_{23}-v_{32}$;

$3^\circ$ $2v_{11}+v_{23}-v_{32}$.

Take an orthonormal basis $\{f_i\}_{1\le i\le 4}$ with $f_1,f_2\in \R^2_1,f_3,f_4\in \R^2_2$, such that
\be \begin{split}w_1=f_1,w_2=-\frac12 f_1+\frac{\sqrt 3}{2}f_2, w_3=-\frac12 f_1-\frac{\sqrt 3}{2}f_2,\\
u_1=f_3,u_2=-\frac12 f_3+\frac{\sqrt 3}{2}f_4, u_3=-\frac12 f_3-\frac{\sqrt 3}{2}f_4.\end{split}\ee

Then for $1^\circ$:
\be  \begin{split}&||2 v_{11}+v_{23}+v_{32}||\\
=&||2f_1\wedge f_3+(-\frac12 f_1+\frac{\sqrt 3}{2}f_2)\wedge(-\frac12 f_3-\frac{\sqrt 3}{2}f_4)+(-\frac12 f_1-\frac{\sqrt 3}{2}f_2)\wedge(-\frac12 f_3+\frac{\sqrt 3}{2}f_4)||\\
=&||\frac52 f_1\wedge f_3-\frac32 f_2\wedge f_4||\le ||f_1||+\frac 32|| f_1\wedge f_3-f_2\wedge f_4||\le 1+\frac 32=\frac 52<3,\end{split}\ee
where the third last inequality is by Lemma 4.1;

For $2^\circ$:\be  \begin{split}&||2 v_{11}-v_{23}-v_{32}||\\
=&||2f_1\wedge f_3-(-\frac12 f_1+\frac{\sqrt 3}{2}f_2)\wedge(-\frac12 f_3-\frac{\sqrt 3}{2}f_4)-(-\frac12 f_1-\frac{\sqrt 3}{2}f_2)\wedge(-\frac12 f_3+\frac{\sqrt 3}{2}f_4)||\\
=&||\frac 32 f_1\wedge f_3+\frac 32 f_2\wedge f_4||=\frac 32||f_1\wedge f_3+f_2\wedge f_4||\le \frac 32<3,\end{split}\ee
where the second last inequality is from Lemma 4.1;

For $3^\circ$: 
\be\begin{split}&||2 v_{11}+v_{23}-v_{32}||\\
=&||2f_1\wedge f_3+(-\frac12 f_1+\frac{\sqrt 3}{2}f_2)\wedge(-\frac12 f_3-\frac{\sqrt 3}{2}f_4)-(-\frac12 f_1-\frac{\sqrt 3}{2}f_2)\wedge(-\frac12 f_3+\frac{\sqrt 3}{2}f_4)||\\
=&||2f_1\wedge f_3+\frac{\sqrt 3}{2}f_1\wedge f_4-\frac{\sqrt 3}{2}f_2\wedge f_3||\\
\le &2||f_1\wedge f_3||+\frac{\sqrt 3}{2}||f_1\wedge f_4-f_2\wedge f_3||\\
\le& 2+\frac{\sqrt 3}{2}<3.
\end{split}\ee

This ends the discuss for Subsubcase 2.1.2, and hence Subcase 2.1.

Subcase 2.2. $\e(2,1)=-\e(3,1).$ Modulo symmetry of indices, we can suppose that $\e(2,1)=\e(1,2)$, and hence $-\e(1,2)v_{11}+\e(2,1)v_{21}+\e(3,1)v_{31}=2\e(2,1)v_{21}$, since $v_{11}+v_{21}+v_{31}=0$. Therefore the sum becomes
\be \sum_{1\le j,l\le 3, l\ne j}\e(j,l)v_{jl}=2\e(2,1)v_{21}+\e(2,3)v_{23}+\e(3,2)v_{32}.\ee

Modulo multiplication of $\pm$, we have also three possibilities:

$1^\circ$ $2v_{21}+v_{23}\pm v_{32}$;

$2^\circ$ $2v_{21}-v_{23}-v_{32}$;

$3^\circ$ $2v_{21}-v_{23}+v_{32}$.

For $1^\circ$, notice that $v_{21}+v_{23}=-v_{22}$, hence the sum becomes $v_{21}-v_{22}\pm v_{32}$, and its norm
\be ||v_{21}-v_{22}\pm v_{32}||\le ||v_{21}||+||v_{22}||+||v_{32}||\le 3.\ee

For the rest two possibilities, first take an orthonormal basis $\{f_i\}_{1\le i\le 4}$ with $f_1,f_2\in \R^2_1,f_3,f_4\in\R^2_2$, such that
\be \begin{split}w_2=f_1,w_3=-\frac12 f_1+\frac{\sqrt 3}{2}f_2, w_1=-\frac12 f_1-\frac{\sqrt 3}{2}f_2,\\
u_1=f_3,u_2=-\frac12 f_3+\frac{\sqrt 3}{2}f_4, u_3=-\frac12 f_3-\frac{\sqrt 3}{2}f_4.\end{split}\ee

Then for $2^\circ$, 
\be \begin{split}||2v_{21}-v_{23}-v_{32}||&=||2f_1\wedge f_3-f_1\wedge(-\frac12 f_3-\frac{\sqrt 3}{2}f_4)-(-\frac12 f_1+\frac{\sqrt 3}{2}f_2)\wedge(-\frac12 f_3+\frac{\sqrt 3}{2}f_4)||\\
&=||\frac94f_1\wedge f_3+\frac{3\sqrt 3}{4}f_1\wedge f_4+\frac{\sqrt 3}{4}f_2\wedge f_3-\frac 34f_2\wedge f_4||\\
&=||\frac{3\sqrt 3}{2}f_1\wedge(\frac{\sqrt 3}{2}f_3+\frac 12f_4)+\frac{\sqrt 3}{2}f_2\wedge(\frac 12 f_3-\frac{\sqrt 3}{2}f_4)||\\
&\le ||\sqrt 3f_1\wedge(\frac{\sqrt 3}{2}f_3+\frac 12f_4)||+||\frac{\sqrt 3}{2}[f_1\wedge(\frac{\sqrt 3}{2}f_3+\frac 12f_4)+f_2\wedge(\frac 12 f_3-\frac{\sqrt 3}{2}f_4)]||.\end{split}\ee
The two vectors $\frac{\sqrt 3}{2}f_3+\frac 12f_4$ and $\frac 12 f_3-\frac{\sqrt 3}{2}f_4$ form an orthonormal basis of $\R^2_2$, hence the four vectors $f_1,f_2, \frac{\sqrt 3}{2}f_3+\frac 12f_4,\frac 12 f_3-\frac{\sqrt 3}{2}f_4$ form an orthonormal basis of $\R^4$. By Lemma 4,1, $||f_1\wedge(\frac{\sqrt 3}{2}f_3+\frac 12f_4)+f_2\wedge(\frac 12 f_3-\frac{\sqrt 3}{2}f_4)||\le 1$, hence by (4.36)
\be ||2v_{21}-v_{23}-v_{32}||\le \sqrt||f_1\wedge(\frac{\sqrt 3}{2}f_3+\frac 12f_4)||+\frac{\sqrt 3}{2}=\sqrt 3+\frac{\sqrt 3}{2}<3.\ee

For $3^\circ$, 
\be \begin{split}||2v_{21}-v_{23}+v_{32}||&=||2f_1\wedge f_3-f_1\wedge(-\frac12 f_3-\frac{\sqrt 3}{2}f_4)+(-\frac12 f_1+\frac{\sqrt 3}{2}f_2)\wedge(-\frac12 f_3+\frac{\sqrt 3}{2}f_4)||\\
&=||\frac{11}{4}f_1\wedge f_3+\frac{\sqrt 3}{4}f_1\wedge f_4-\frac{\sqrt 3}{4}f_2\wedge f_3+\frac 34 f_2\wedge f_4||\\
&=||2 f_1\wedge f_3+\frac{\sqrt 3}{2}[f_1\wedge (\frac{\sqrt 3}{2}f_3+\frac 12 f_4)+f_2\wedge(- \frac12 f_3+\frac{\sqrt 3}{2} f_4)]||\\
&\le 2||f_1\wedge f_3||+\frac{\sqrt 3}{2}||f_1\wedge (\frac{\sqrt 3}{2}f_3+\frac 12 f_4)+f_2\wedge(- \frac12 f_3+\frac{\sqrt 3}{2} f_4)||\\
&=2+\frac{\sqrt 3}{2}||f_1\wedge (\frac{\sqrt 3}{2}f_3+\frac 12 f_4)+f_2\wedge(- \frac12 f_3+\frac{\sqrt 3}{2} f_4)||.
\end{split}\ee
The two vectors $\frac{\sqrt 3}{2}f_3+\frac 12 f_4$ and $- \frac12 f_3+\frac{\sqrt 3}{2} f_4$ form an orthonormal basis of $\R^2_2$, hence $\{f_1,f_2,\frac{\sqrt 3}{2}f_3+\frac 12 f_4,- \frac12 f_3+\frac{\sqrt 3}{2} f_4$ form an orthonormal basis of $\R^4$. By Lemma 4,1, $||f_1\wedge (\frac{\sqrt 3}{2}f_3+\frac 12 f_4)+f_2\wedge(- \frac12 f_3+\frac{\sqrt 3}{2} f_4)||\le 1$, hence
\be ||2v_{21}-v_{23}+v_{32}||\le 2+\frac{\sqrt 3}{2}<3.\ee

We have discussed all possible cases, and hence the proof for Lemma 3.24 is completed. \qed
 
 \section{Almgren and topological minimality of $Y\times Y$}
 
 In this section we will establish the Almgren minimality of $Y\times Y$, and will prove some sort of topological minimality as well. 
 
 \subsection{Almgren minimality of $Y\times Y$}
 
 Recall that the definition of Almgren minimal sets is given in Definition 0.1.
 
 \begin{thm}\label{al}The set $Y\times Y$ is an Almgren minimal cone in $\R^4$.
 \end{thm}
 
 \nd Since $Y\times Y$ is a cone, it is Almgren minimal if and only if $Y\times Y$ minimizes the Hausdorff measure among all its Almgren competitors in $D^\circ$. Hence it suffices to prove that any local Almgren competitor of $Y\times Y$ in $D^\circ$ belongs to $\F$.
 
 But this is natural: for any Almgren competitor $F$ of $Y\times Y$ in $D^\circ$, there exists $f:\R^4\to\R^4$, such that $f=Id$ outside some compact set $K\subset D^\circ$, and $f(K)\subset K$. In particular, $f$ induces a homomorphism between the singular homology groups $f_*: H_1(Y\times Y\cap D)\to H_1(F\cap D)$. Since $f$ fixes all the $\gamma_{jl}$, $\gamma_{jl}$ represents a zero element in $H_1(F\cap D, \Z/2\Z)$, hence $F\in \F$.
 
This yields that $Y\times Y$ is an Almgren minimal set in $\R^4$. \qed

\subsection{Topological minimality of $Y\times Y$}

Next we deal with the topological minimality of $Y\times Y$. The essential point of the proof is the Poincar\'e-Lefschetz duality for homology and cohomology groups.

Let us first give the definition of topological minimal sets. The notion of topological minimality that we use here is almost the same as that was introduced in \cite{topo}, except that here we use $\Z/2\Z$ as the coefficient group, rather than $\Z$. (Hence in the rest of the article, if not pointed out, the coefficient group is automatically $\Z/2\Z$). However, all the properties that were proved in \cite{topo} can also be proved in the same way when we replace $\Z$ by $\Z/2\Z$. 

\begin{defn}[Topological competitors (with coefficient in $\Z/2\Z$)] Let $E$ be a closed set in $\R^n$. We say that a closed set $F$ is a topological competitor of dimension $d$ ($d<n$) of $E$ (with coefficient in $\Z/2\Z$), if there exists a ball $B\subset\R^n$ such that

1) $F\bs B=E\bs B$;

2) For each finite sum of disjoint Euclidean $n-d-1$-sphere $S=S_1+S_2+\cdots+S_k\subset\R^n\bs(B\cup E)$, if $S$ represents a non-zero element in the singular homology group $H_{n-d-1}(\R^n\bs E;\Z/2\Z)$, then it is also non-zero in $H_{n-d-1}(\R^n\bs F;\Z/2\Z)$.

We also say that such a $F$ is a topological competitor of $E$ with respect to the ball $B$.
\end{defn}

\begin{defn}[Topological minimal sets (with coefficient in $\Z/2\Z$)]\label{min}
Let $0<d<n$ be integers. A closed set $E$ in $\R^n$ is said to be topologically minimal (with coefficient in $\Z/2\Z$) of dimension $d$ in $\R^n$ if 
\be H^d(E\cap B)<\infty\mbox{ for every compact ball }B\subset U,\ee
and
\be H^d(E\bs F)\le H^d(F\bs E)\ee
for every topological competitor $F$ of dimension $d$ (with coefficient in $\Z/2\Z$) for $E$.
\end{defn}

\begin{rem}The definition is the same if we replace "singular homology" by "simplicial homology", because $\R^n\bs E$ and $\R^n\bs F$ are $C^\infty$ $n-$dimensional manifolds, which admit locally finite triangulations, and hence the two homology groups are naturally isomorphic. (cf.\cite{Hat} Theorem 2.27).
\end{rem}

\begin{rem}The exact same proof of Corollary 3.17 in \cite{topo} gives that all topological minimal sets with coefficients in $\Z/2\Z$ are Almgren minimal sets.
\end{rem}

 \begin{thm}The set $Y\times Y$ is a topological minimal cone (with coefficient in $\Z/2\Z$) in $\R^4$.
 \end{thm}
 
 \nd Denote by $\F_T$ the class of all topological competitors $F$ of $Y\times Y$ with respect to the ball $B=B(0,\frac12)\subset\R^4$. Set
 \be \F_T^R=\{F\in\F_T; F\mbox{ is regular}\}.\ee
 Then by the same argument as in Proposition 1.4, and notice that the every deformation (in particular a Federer Fleming projection) preserves the topological condition 2) in Definition 5.2 (cf. \cite{topo} Proposition 3.7), we have (notice that $B\subset D$)
 \be \inf\{H^2(F\cap D); F\in \F_T\}=\inf\{H^2(F\cap D); F\in \F_T^R\}.\ee
 
 So in order to prove Theorem 5.8, it suffices to show that 
 \be\F_T^R\subset\F^R.\ee
 
 Let $F\in\F_T^R$ be a regular topological competitor of $Y\times Y$ in the ball $B$. Then by definition of topological competitors, and the fact that $B\subset D$, we have $F\bs D^\circ=Y\times Y\bs D^\circ$. So in order to get (5.11), we only have to check that for $1\le j,l\le 3$, $\gamma_{jl}$ represents a zero element in the simplicial homology group $H_1(F\cap D,\Z/2\Z)$.
 
The set $Y\times Y$ is the union of nine $\frac14-$planes $H_{ik}:=[o_1,a_i)\times[o_2,b_k), 1\le i,k,\le 3$. Set $c_{ik}=(a_i,b_k)\in\R^4$, the point whose first two coordinates coincide with the coordinates of $a_i$, and the last two coordinates coincide with those of $b_k$. Denote by $\Pi_{ik}$ the plane in $\R^4$ that is orthogonal to $H_{ik}$, and passing through $c_{ik}$. Denote by $s_{ik}$ the circle of radius $\frac{1}{10}$, centered at the point $c_{ik}$, and that lies in $\Pi_{ik}$. Then $s_{ik}$ is a small circle outside $D$ that links $H_{ik}$.We denote also by $s_{ik}$ the corresponding element in homology groups with coefficients in $\Z/2\Z$ for short. Then  in $H_1(\R^4\bs F,\Z/2\Z)$
 \be s_{i1}+s_{i_2}+s_{i3}=0=s_{1k}+s_{2k}+s_{3k}, \mbox{ for }1\le i,k\le 3,\ee
 and, in the special case of $Y\times Y$, $H_1(\R^4\bs Y\times Y,\Z/2\Z)$ is the Abelian group generated by $s_{ik},1\le i,k\le 3$ with the relations (5.12). Notice that the relations (5.12) has in fact only 5 independent relations. Thus $H_1(\R^4\bs Y\times Y,\Z/2\Z)$ is in fact a vector space (since $\Z/2\Z$ is a field) with basis $\{s_{ik},1\le i,k\le 2\}$.
 
 Take $j=1,l=1$ for example. We want to show that $\gamma_{11}$ is zero in $H_1(F\cap D,\Z/2\Z)$.
 
 Denote by $V=[o_1,a_1)\cup[o_1,a_2)$ the union of two branches of $ Y_1$, and $U=[o_2,b_1)\cup[o_2,b_2)\subset Y_2$ the union of two branches of $Y_2$. Hence $V\cap B_1=S_1,U\cap B_2=R_1$, and $V\times U\subset\R^4$ is topologically a plane. Denote $E=(F\cap D^\circ)\cup(V\times U\bs D^\circ)$. It is easy to see that the four circles $s_{11},s_{12},s_{21}$ and $s_{22}$ represent a same element in $H_1(\R^4\bs E,\Z/2\Z)$. Denote by $s$ this element in $H_1(\R^4\bs E,\Z/2\Z)$.
 
 We want to show that $s\ne 0$. Suppose not, that is, $s=0$ in $H_1(\R^4\bs E,\Z/2\Z)$. Then there exists a $C^1$ simplicial 2-chain $\Gamma$ in $\R^4\bs E$ such that $\partial \Gamma=s_{11}$. Since $F$ is a topological competitor of $Y\times Y$ in $B$, $s_{11}\ne 0$ in $H_1(\R^4\bs F,\Z/2\Z)$, hence $\Gamma\cap F\ne\emptyset$. But $\Gamma\subset\R^4\bs E$, hence $\Gamma$ can only meet F at $[H_{13}\cup H_{23}\cup H_{31}\cup H_{32}\cup H_{33}]\bs D$. We can also ask that $\Gamma$ meet these five $\frac14$ planes transversally, and do not meet any of their intersections. This yields that $\partial\Gamma$ is a sum 
 \be\partial\Gamma=s_{11}+\d_{13}s_{13}+\d_{23}s_{23}+\d_{33}s_{33}+\d_{31}s_{31}+\d_{32}s_{32},\ee
 with $\d_{ik}\in\{0,\pm 1\}$. Moreover, at least one of these five $\d_{ik}$ is non zero. This gives that in $H_1(\R^4\bs F,\Z/2\Z)$,
 \be s_{11}+\d_{13}s_{13}+\d_{23}s_{23}+\d_{33}s_{33}+\d_{31}s_{31}+\d_{32}s_{32}=0.\ee
 
Combine with (5.12),we have
 \be s_{11}+\d_{13}[s_{11}+s_{12}]+\d_{23}[s_{21}+s_{22}]+(\d_{31}+\d_{33})[s_{21}+s_{11}]+(\d_{32}+\d_{33})[s_{12}+s_{22}]=0.\ee
 We simplify and get
 \be [1+\d_{13}+\d_{31}+\d_{33}]s_{11}+[\d_{13}+\d_{32}+\d_{33}]s_{12}+[\d_{23}+\d_{31}+\d_{33}]s_{21}+[\d_{23}+\d_{32}+\d_{33}]s_{22}=0,\ee
in $H_1(\R^4\bs F,\Z/2\Z)$.

But $s_{11},s_{12},s_{21},s_{22}$ are independent elements in $H_1(\R^4\bs Y\times Y,\Z/2\Z)$, hence (5.16) gives
\be 1+\d_{13}+\d_{31}+\d_{33}=\d_{13}+\d_{32}+\d_{33}=\d_{23}+\d_{31}+\d_{33}=\d_{23}+\d_{32}+\d_{33}=0.\ee
However, the sum of the four numbers gives
\be \begin{split}0&=[1+\d_{13}+\d_{31}+\d_{33}]+[\d_{13}+\d_{32}+\d_{33}]+[\d_{23}+\d_{31}+\d_{33}]+[\d_{23}+\d_{32}+\d_{33}]\\
&=1+2(\d_{13}+\d_{31}+\d_{32}+\d_{23}+2\d_{33})=1,
 \end{split}\ee
 which is impossible. Hence we have a contradiction.
 
 Hence in $H_1(\R^4\bs E,\Z/2\Z)$, $s_{11}=s\ne 0$ in $H_1(\R^4\bs E,\Z/2\Z)$.
 
 Now since $H_1(\R^4\bs E,\Z/2\Z)$ is a finite dimensional vector space (because $\Z/2\Z$ is a field), $H^1(\R^4\bs E,\Z/2\Z)=[H_1(\R^4\bs E,\Z/2\Z)]^*$. For the non zero element $s_{11}\in H_1(\R^4\bs E,\Z/2\Z)$, denote by $\xi$ the dual element of $s_{11}$ in the cohomology group $H^1(\R^4\bs E,\Z/2\Z)$. Denote by $\varphi$ the isomorphism of the Poincar\'e-Lefschetz duality
 \be\varphi:H^1(\R^4\bs E,\Z/2\Z)\cong H_3(\R^4,E,\Z/2\Z),\ee
 Then $\varphi(s_{11})$ can be represented by a smooth simplicial 3-chain $\Sigma$ with boundary in $E$, and $|\Sigma|\cap s_{11}$ is a single point. Denote by $\eta=\partial\Sigma$, then this is a 2-chain with support in $E$ such that $s_{11}$ is non-zero in $H_1(\R^4\bs |\eta|,\Z/2\Z)$. 
 
 Notice that outside $D$, the set $E$ is topologically a plane, which is linked by $s_{11}$, hence if $s_{11}$ is non-zero in $H_1(\R^4\bs |\eta|,\Z/2\Z)$, then $|\eta|\supset (E\bs D^\circ)$. 
  
 But $E\cap\partial D=\gamma_{11}$, hence $\gamma_{11}=\partial(E\bs D^\circ)=\partial (\eta-(E\bs D^\circ))$. Notice that the support of $\eta-(E\bs D^\circ)$ is contained in $E\cap D$, hence $\gamma_{11}$ is a boundary in $E\cap D$, which yields that $\gamma_{11}$ represents a zero element in $H_1(E\cap D,\Z/2\Z)$, and thus in $H_1(F\cap D,\Z/2\Z)$.
 
 The same arguments holds for all $\gamma_{jl},1\le j,l\le 3$. Thus we have proved that $F\in \F^R$. Since this is true for any regular topological competitor $F$ of $Y\times Y$, we have (5.11), and thus the proof of Theorem 5.8 is completed.\qed
 
 \begin{cor}For any $k>0$, the product $Y\times Y\times \R^k$ is a $2+k$-dimensional topological and Almgren minimal in $\R^{k+4}$.
 \end{cor}
 \nd The set $Y\times Y\times \R^k$ is the product of a topological minimal set $Y\times Y$ with $\R^k$, hence it is topologically minimal. (cf., the proof of Proposition 3.23 of \cite{topo}). But the topological minimality implies the Almgren minimality (cf. the proofs of Proposition 3.7 and Corollary 3.17 of \cite{topo}, hence the set is also Almgren minimal.\qed
 
 \section{Minimality for $Y\times Y$ is not stable by small perturbation}
 
 In this section we want to prove that there is no other minimal cone in $\R^4$ that admits the same topology as $Y\times Y$. In particular, small perturbations around $Y\times Y$ do not preserve the minimality.
 
 Let us first recall a necessary condition for 2-dimensional cones to be Almgren-minimal in $\R^n$. Let $C\in\R^n$ be a 2-dimensional Almgren minimal cone. Then by \cite{DJT} Proposition 14.1, 
 \be\begin{split}&\mbox{The intersection of }C\mbox{ and the unit sphere }X=C\cap\partial B(0,1)
 \mbox { is a finite union }\\
 &\mbox{of great circles or arcs of great circles.}\mbox{ The circles are far from the rest of }X.\\
 &\mbox{ At their ends, the arcs meet by set of three, with }120^\circ\mbox{ angles (no free ends).} \end{split}\ee
 
 \begin{rem}It is useful to know that any three arcs of great circles in $\mathbb S^3$ that meet at a common point with $120^\circ$ angle belong to a same (2-dimensional) great sphere, i.e. the intersection of $\mathbb S^3$ with a 3-dimensional vector space.
 \end{rem}
 
 Our cone $Y\times Y$ satisfies obviously the condition (6.1). And we are going to prove the following theorem.
 
 \begin{thm} The only cone (modulo isometry) with the same topology as $Y\times Y$ in $\R^4$ that satisfies the condition (6.1) is $Y\times Y$.\end{thm}
 
 \nd
 
 Let us first look at the structure of $K=Y\times Y\cap \partial B(0,1)$. The set K is a net that is composed of 9 arcs of great circles $K_{ij}, 1\le i,j\le 3$, and has six $\Y$ points $x_1,x_2,x_3,y_1,y_2,y_3$. The endpoints of $K_{ij}$ are $x_i$ and $y_j$. Modulo a rotation we can suppose that $y_3=(0,0,0,1)$.
 
 Let us just work on the unit sphere $\mathbb S^3$ of $\R^4$. We take the stereographic projection $\pi$, with projection point $(0,0,0,1)$, from $\mathbb S^3\bs \{(0,0,0,1)\}$ to the 3-plane $\pi=\{(x,y,z,t): t=-1\}$, which is the plane tangent to the unit sphere at $(0,0,0,-1)$. The map $\pi$ sends the point $y_3$ to the $\infty$ point. Then $\pi(K)$ is, topologically, like in Figure 2 below.
 
 \centerline{\includegraphics[width=1.0\textwidth]{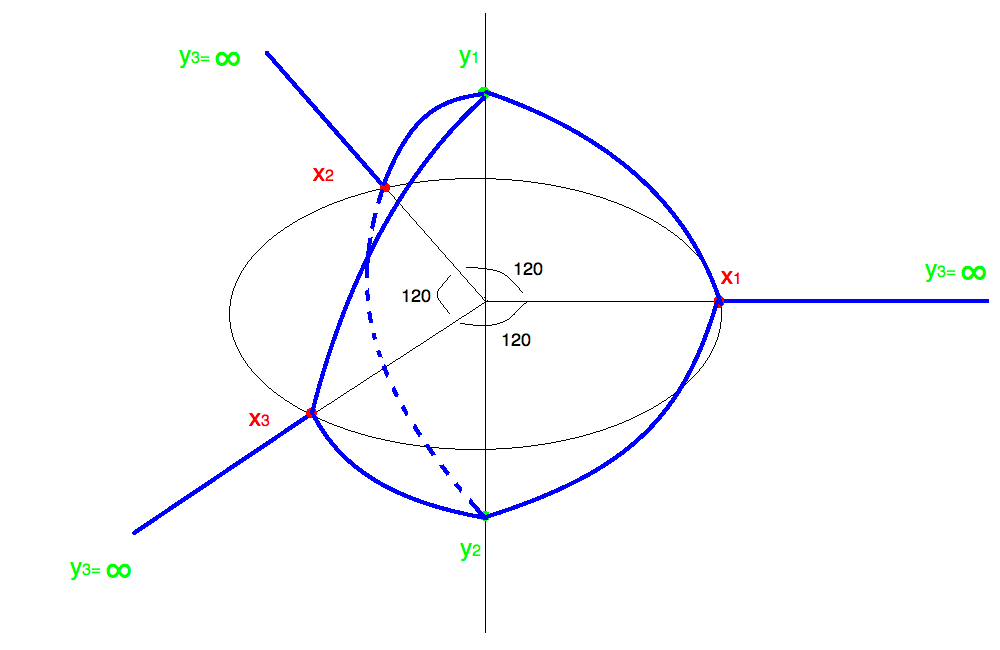}}
\nopagebreak[4]
\centerline{Figure 2} 
 
 Now let $C$ be another cone in $\R^4$ that admits the same topology as $Y\times Y$ and satisfies (6.1). Then $X=C\cap\partial B(0,1)$ is composed of 9 arcs of great circles $X_{ij},1\le i,j\le 3$, and has six $\Y$ points $a_1,a_2,a_3,b_1,b_2,b_3$. The endpoints of $X_{ij}$ is $a_i$ and $b_j$, $1\le i,j\le 3$. Without loss of generality, suppose that $b_3$ is the point $(0,0,0,1)$. We take also the same stereographic projection $\pi$ as above.
 
 Now since the arcs $X_{i3},1\le i\le 3$ meet at the point $b_3$ with $120^\circ$ angles, by Remark 6.2, the four points $a_1,a_2,a_3,b_3$ belong to a same great sphere $S$, and if we prolong the arcs $X_{i3}$, they will meet at the opposite point $(0,0,0,-1)$, whose image by $\pi$ is the origin. Since the image of any great sphere passing through the point $\pi(b_3)=\infty$ is a plane in $\R^3$ passing through the origin, $\pi(a_1),\pi(a_2),\pi(a_3)$ belong respectively to the three branches of a 1-dimensional $\Y$ set $Y$ in $\R^3$ that is centered at the origin. See Figure 3. Denote by $Q=\pi(S)$ the plane containing $Y$. 
 
  \centerline{\includegraphics[width=1.0\textwidth]{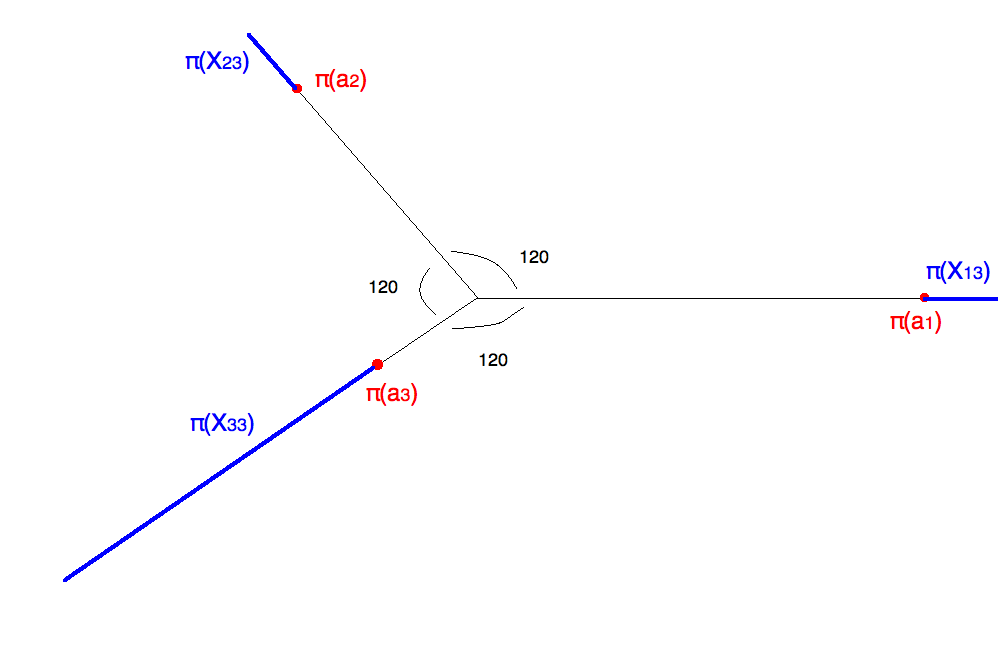}}
\nopagebreak[4]
\centerline{Figure 3} 

Now the same is true for $b_1$ and $b_2$, that is, $\{a_i\}_{1\le i\le 3}$ and $b_1$ belong to a same great sphere, and $\{a_i\}_{1\le i\le 3}$ and $b_2$ belong also to a great sphere.

We claim that
\be \mbox{the three points }\{a_i\}_{1\le i\le 3}\mbox{ belong to a same great circle}.\ee

In fact if they do not belong to a same great circle, then there is a unique great sphere $S'$ that contains them three. In this case, $b_1,b_2$ belong to $S'$, too. But we already have a great sphere $S$ that contain $\{a_i\}_{1\le i\le 3}$ and $b_3$, hence all the six points $a_i,b_j, 1\le i,j\le 3$ belong to the 2-dimensional sphere $S$. But this is impossible, because topologically, the net $X$ is the graph $K_{3,3}$, which is not a planar graph.

Thus we have the claim (6.4).

It is easy to check that any two of the three $a_i,1\le i\le 3$ are not opposite points in $\R^4$, since the images by $\pi$ of the three points belong to three different branches of a 1-dimensional $\Y$ set, and the images by $\pi$ of  any two opposite points in $\mathbb S^3$ should belong to a line that passing through the origin.

Hence in fact, any two of them define already a unique great circle, and hence the position of any two of them fixes also the third.

Now for $i=1,2,3$, since the arcs $X_{ij},1\le j\le 3$ meet at the point $a_i$ with $120^\circ$, by Remark 6.2, $a_i,b_1,b_2,b_3$ belong to a same great sphere. But the image by $\pi$ of any great sphere that contains the point $b_3$ is a plane passing through the origin, hence the points $\pi(a_i),\pi(b_1),\pi(b_2)$ and the origin $o$ belong to a same plane $P_i\subset\R^3$.  Denote by $D\subset\R^3$ the line that contains $\pi(b_1)$ and $\pi(b_2)$, then $P_1\cap P_2\cap P_3\supset D$. On the other hand, since $\pi(a_1)\in P_1\bs P_2$, therefore $P_1\ne P_2$, which yields, $P_1\cap P_2$ is a line, and hence $P_1\cap P_2\cap P_3\subset D$. As a result, $P_1\cap P_2\cap P_3=D$.  Hence $D$ contains the origin. 
%(Or, we can get, by the same argument as $\{a_i\}_{1\le i\le 3}$, that the three points $\{b_i\}_{1\le i\le 3}$ belong to a same great circle, whose image by $\pi$ is a line that passing through the origin.) 
Denote by $P_i'$ the half plane in $P_i$ bounded by $D$ and containing the half line $\pi(X_{i3})$. We have also $\cap_{1\le i\le 3}P_i'=D$.

Notice that $\pi$ is a conformal mapping, hence the curves $\pi(X_{11}),\pi(X_{21}),\pi(X_{31})$ meet at the point $b_1$ with $120^\circ$. This means, if we denote by $L_i,1\le i\le 3$ the half line in $P_i'$ that is bounded by $b_1$ and tangent to the curve $\pi(X_{i1})$ at the point $b_1$, then the three half lines meet at $b_1$ with angles of $120^\circ$. (That is, $\cup_{1\le i\le 3}L_i$ is a one dimensional $\Y$ set). On the other hand, $X_{11},X_{21},X_{31}$ belong to a same great sphere, which means their images $\pi(X_{11}),\pi(X_{21}),\pi(X_{31})$ also belong to a sphere $S_1$ in $\R^3$, hence the tangent half lines $L_i,1\le i\le 3$ belong to the  plane $Q_1$ that is tangent to $S_1$ at $b_1$.  

The same argument gives that there is a plane $Q_2$ containing $b_2$, such that the intersection lines $L_i',1\le i\le 3$ of $Q_2$ with $P_i'$ makes $120^\circ$ angles (i.e. $\cup_{1\le i\le 3}L_i'$ is a one dimensional $\Y$ set). $Q_2$ is also the tangent plane at $b_1$ to the $\pi-$image of the great sphere that contains $X_{12},X_{22},X_{32}$.

Let us sum up a little: there are three half planes $P_i',1\le i\le 3$, which are contained in three different planes $P_i,1\le i\le 3$, who meet along a line $D$. There are two planes $Q_1,Q_2$, $Q_1\cap P_i=L_i,Q_2\cap P_i=L_i'$. The three $L_i,1\le i\le 3$ make $120^\circ$ angles between each other, and so do $L_i',1\le i\le 3$.

We claim that, in $\R^3$
\be Q_1\mbox{ and }Q_2\mbox{ are, either parallel, either symmetric with respect to the plane }D^\perp,\ee  
where $D^\perp$ is the plane that is orthogonal to the line $D$.

Denote by $w$ a unit direction vector of $D$. Denote by $v_i,1\le i\le 3$ the unit vector such that $v_i+b_1\in P_i'$ and $v_i\perp w$. See Figure 4.

  \centerline{\includegraphics[width=0.5\textwidth]{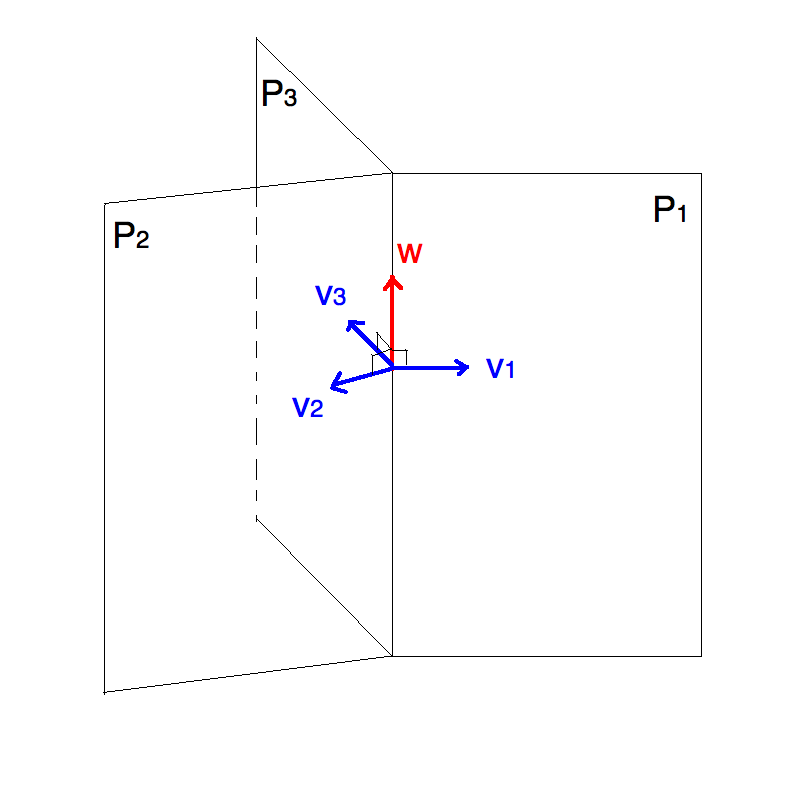}}
\nopagebreak[4]
\centerline{Figure 4} 

Now for $Q_1$, take the three unit vectors $q_i, 1\le i\le 3$ such that $L_i=\{b_1+ rq_i,r\ge 0\}$ (i.e, $q_i$ is the direction vector of $L_i$. Then $b_1+q_i$ belong to the half plane $P_i'$. Hence there exists $\theta_i\in (-\frac\pi2,\frac\pi2),1\le i\le 3$ such that
\be q_i=\cos\theta_i v_i+\sin\theta_i w,1\le i\le 3.\ee

The three $q_i,1\le i\le 3$ makes $120^\circ$ angles, hence their sum is zero. Thus

\be 0=\sum_{1\le i\le 3}q_i=\sum_{1\le i\le 3}\cos\theta_i v_i+\sin\theta_i w=\sum_{1\le i\le 3}\cos\theta_i v_i+[\sum_{1\le i\le 3}\sin\theta_i ]w.\ee

Notice that $w$ is independent with the plane generated by $v_i,1\le i\le 3$, hence
\be \sum_{1\le i\le 3}\sin\theta_i=0,\ee
and
\be\sum_{1\le i\le 3}\cos\theta_i v_i=0.\ee

Now since any two of the three $v_i,1\le i\le 3$ are independent, there is only one possible pair of real numbers $(t_2,t_3)$ (which, of course, depends on $w,v_i,1\le i\le 3$, and hence $P'_i,1\le i\le 3$) such that
\be v_1+t_2v_2+t_3v_3=0.\ee
And, modulo a possible permutation of indices, we can suppose that $1\le |t_2|\le |t_3|.$
Hence by (6.9), there exists a unique constant $T\in[0,1]$ such that
\be \cos\theta_1=T,\cos\theta_2=Tt_2,\cos\theta_3=Tt_3.\ee

Combine this with (6.8), we get that there exists $\e_i\in\{\pm 1\}$ such that
\be \sqrt{1-T^2}+\e_2\sqrt{1-T^2t_2^2}+\e_3\sqrt{1-T^2t_3^2}=0.\ee

Now we know that $1\le t_2^2\le t_3^2$, hence $\sqrt{1-T^2}\ge \sqrt{1-T^2t_2^2}\ge\sqrt{1-T^2t_3^2}\ge 0$, which means $\e_2=\e_3=-1$. Hence in fact $\theta_i,1\le i\le 3$ satisfy
\be \sqrt{1-T^2}-\sqrt{1-T^2t_2^2}-\sqrt{1-T^2t_3^2}=0.\ee

 Define $f(u)=\sqrt{1-u}-\sqrt{1-ut_2^2}-\sqrt{1-ut_3^2}$. Then $T^2$ is a solution for $f(u)=0$ in the interval $[0,1]$.
 
 But 
 \be \begin{split}f'(u)&=\frac12[\frac{-1}{\sqrt{1-u}}+\frac{t_2^2}{\sqrt{1-ut_2^2}}+\frac{t_3^2}{\sqrt{1-ut_3^2}}]\\
 &=\frac12\left[\frac{-\sqrt{1-ut_2^2}\sqrt{1-ut_3^2}+t_2^2\sqrt{1-u}\sqrt{1-ut_3^2}+t_3^2\sqrt{1-u}\sqrt{1-ut_2^2}}{\sqrt{1-u}\sqrt{1-ut_2^2}\sqrt{1-ut_3^2}}\right],\end{split}\ee
it is always strictly larger than 0. Hence the equation $f(u)=0$ admits only one solution $u_0$, which means that $T=\sqrt{u_0}$, since $T>0$.

Then for $\theta_i,1\le i\le 3$ we have two solutions :
\be
\left\{\begin{array}{ccc}\theta_1&=&\arcsin\sqrt{1-u_0}, \\\theta_2&=&-\arcsin\sqrt{1-u_0t_2^2}, \\ \theta_3&=&-\arcsin\sqrt{1-u_0t_3^2}\end{array}\right. ,
\mbox{ or }
\left\{\begin{array}{ccc}\theta_1&=&-\arcsin\sqrt{1-u_0}, \\\theta_2&=&\arcsin\sqrt{1-u_0t_2^2}, \\ \theta_3&=&\arcsin\sqrt{1-u_0t_3^2}\end{array}\right. . \ee

The same argument works also for $Q_2$, where we denote by $q_i'$ the three unit direction vectors of $L_i'$, and $\theta_i'$ the angles such that \be q_i'=\cos\theta_i' v_i+\sin\theta_i' w,1\le i\le 3,\ee
Thus, either $\theta_i=\theta_i'$ for all $1\le i\le 3$, either $\theta_i=-\theta_i'$ for all $1\le i\le 3$. Geometrically, the first case means $Q_2$ is parallel to $Q_1$, and the second case means $Q_1$ and $Q_2$ are symmetric with respect to $D^\perp$. Thus we get our claim (6.5).

On the other hand, notice that the arc $X_{11}$ (the arc with endpoints $a_1$ and $b_1$) is part of a great circle $C_1$ of $\mathbb S^3$. The great circle $C_1$ passes through the opposite point $-a_1$ of $a_1$, since it passes through $a_1$. 
%A simple calculation gives 
%\be \pi(-a_1)=\frac{-4}{|\pi(a_1)|^2}\pi(a_1).\ee

But $X_{11}$ also lies in the great sphere that contains $a_1,b_1,b_3$, hence $\pi(X_{11})\subset P_1$, which means $\pi(C_1)\subset P_1$. Hence in fact $\pi(C_1)$ is a planar circle that contains the two points $\pi(a_1),\pi(-a_1)$, and is contained in the plane $P_1$. 
%
%  \centerline{\includegraphics[width=1.0\textwidth]{C1.jpg}}
%\nopagebreak[4]
%\centerline{Figure 5} 

We also know that the vector $\overrightarrow{\pi(-a_2)\pi(a_1)}$ and the curve $\pi(X_{11})$ make a $120^\circ$ angle at the point $\pi(a_1)$. Denote by $L$ the tangent line to $\pi(C_1)$ at the point $\pi(a_1)$. Take a point $y$ on $L$ that lies outside the disc bounded by $\pi(C_2)$. Then $\angle_{y\pi(a_1)\pi(-a_1)}=120^\circ$. 

Denote by $a\in P_1$ the  midpoint of  $\pi(a_1)$ and $\pi(-a_1)$, $Z\subset P_1$ the line passing through $a$ and perpendicular to the segment $[\pi(a_1),\pi(-a_1)]$, and denote by $z_1$ the intersection point of $Z$ with the arc $A_1$ of $\pi(C_1)$ with endpoints $\pi(a_1),\pi(-a_1)$ and containing $\pi(X_{11})$, see Figure 5 below: 

  \centerline{\includegraphics[width=0.6\textwidth]{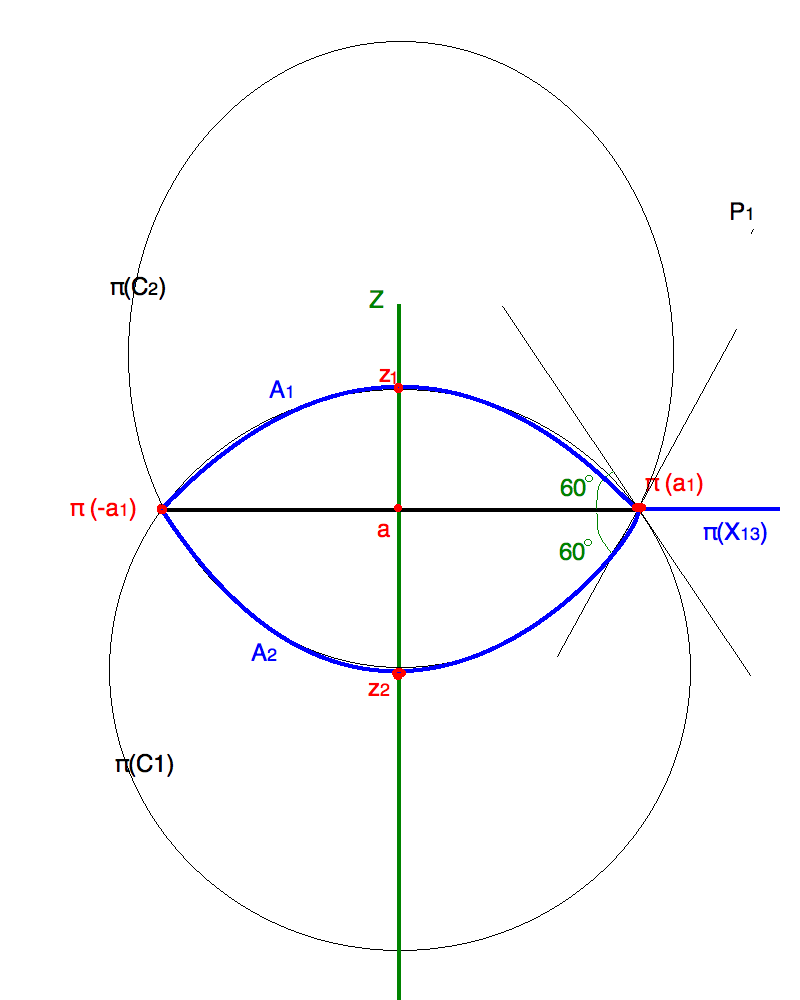}}
\nopagebreak[4]
\centerline{Figure 5}

Then $\angle_{y\pi(a_1)\pi(-a_1)}=120^\circ$ yields that the arc of circle $\wideparen{\pi(a_1)z_1\pi(-a_1)}$ is a $120^\circ$ arc.
Hence the angle $\angle \pi(a_1)z_1\pi(-a_1)=120^\circ$.

We have the same conclusion for the point $b_2$, that is, if we denote by $C_2$ the planar circle in $P_1$ that contains $X_{12}$, $z_2$ the intersection point of $Z$ with the arc  $A_2$ of $\pi(C_2)$ with endpoints $\pi(a_1),\pi(-a_1)$ and containing $\pi(X_{12})$, then $\angle \pi(a_1)z_2\pi(-a_1)=120^\circ$. (Figure 5).

This means, in fact $z_2$ is the center of $C_1$ and $z_1$ is the center of $C_2$, $C_1$ and $C_2$ are symmetric with respect to the segment $[\pi(a_1),\pi(-a_1)]$.

Notice that $b_1\in A_1$, and $b_2\in A_2$. 

Now we know that the two lines $M_1=Q_1\cap P_1$ is tangent to $A_1$ at $b_1$, and $M_2=Q_2\cap P_2$ is tangent to $A_2$ at $b_2$. We have the following two cases:

Case 1: If $Q_1$ and $Q_2$ are parallel, then the two lines $M_1,M_2$ are parallel, too, in this case everything are symmetric, in particular, $a=\frac 12(\pi(a_1)+\pi(-a_1))$ is on the line passing through $b_1,b_2$, i.e. $a\in D$. Notice that $D\cap [\pi(a_1),\pi(-a_1)]=o$, hence $\pi(a_1)$ and $\pi(-a_1)$ are symmetric with respect to the origin. The same argument works also for $a_2,a_3$, which gives
\be o=\frac 12(\pi(a_1)+\pi(-a_1))=\frac 12(\pi(a_2)+\pi(-a_2))=\frac 12(\pi(a_3)+\pi(-a_3)).\ee
This is exactly the case for $Y\times Y$.

Case 2: If $Q_1$ and $Q_2$ are symmetric with respect to $D^\perp$, then since $D\subset P_1$, the two lines $M_1$ and $M_2$ are also symmetric with respect to $D^\perp\cap P_1$. Notice that $D$ is the line passing through $b_1$ and $b_2$, hence the line $D$ has to be perpendicular to the segment $[\pi(a_1),\pi(-a_1)]$. Similarly $D$ has to be perpendicular to the segment $[\pi(a_2),\pi(-a_2)]$, hence $D$ is perpendicular to the plane $Q$. (Recall that $Q$ is the plane containing $\pi(a_i),1\le i\le 3$ and the origin, $Y\subset Q$, and the three $\pi(a_i)$ belong to different branches of $Y$.)

In this case, the three half planes $P_1',P_2',P_3'$ make $120^\circ$ themselves. In this particular case, if we go back to the calculation that was used to prove the claim (6.5), where we have $t_1=t_3=1$, $\theta_1=\theta_2=\theta_3$, hence $T^2=1$, and thus $\sin\theta_i=0,i=1,2,3$. But in this case we have $\theta_i=\theta_i'=0$. This gives again our set $Y\times Y$.

We get therefore $C=Y\times Y$, and the proof of Theorem 6.3 is completed. \qed

 \section{Product of two sets}
 
 As we said in the introduction, we prove the following related theorem in this section.
 
 \begin{pro}Let $E_i,i=1,2$ be two closed rectifiable set of dimension $d_i$ in $\R^{n_i}$ respectively. Set $n=n_1+n_2,d=d_1+d_2$. Then if $E=E_1\times E_2\subset \R^n$ is a $d$-dimensional Almgren minimal set in $\R^n$, then $E_i$ is Almgren minimal of dimension $d_i$ in $\R^{n_i}$, $i=1,2$.
\end{pro}

\nd In all that follows, minimal set means Almgren minimal set.

First we can suppose that both $E_1$ and $E_2$ are reduced, that is, for $i=1,2$, for any point $x\in E_i$ and any $r>0$, $H^{d_i}(B(x,r)\cap E_i)>0$. Otherwise we can replace $E_i$ by its closed support $E_i^*$ (the reduced set $E_i^*\subset E$ with $H^2(E\bs E^*)=0$). Notice that a closed set is minimal if and only if its closed support is minimal.

We will follow the argument in \cite{DJT} proposition 8.3, with some modifications.

So let $f$ be any deformation of $E_1$ in a ball $B\subset \R^{n_1}$  (See Definition 0.1). Set
\be c=H^{d_1}(E_1\cap B)-H^{d_1}(F\cap B).\ee
We want to show that $c\le 0$. For this purpose, we will construct a Lipschitz deformation $\varphi$ from $\R^n$ to $\R^n$, and use the fact that $E$ is minimal.

Let $R>0$ be large, to be decided later. Take a smooth map $\psi$ on $\R^{n_2}$, with
\be\begin{array}{rcl}
    \psi(y)=1&\mbox{ for }&|y|\le R;\\
0\le\psi(y)\le 1&\mbox{ for }&R<|y|<R+1;\\
\psi(y)=0&\mbox{ for }&|y|\ge R+1,
   \end{array}\ee
and $|\nabla\psi(y)|\le 2$ everywhere. Next we define $g:\R^n\to\R^{n_1}$ as follows:
\be g(x,y)=\psi(y)f(x)+(1-\psi(y))(f(x)-x)\mbox{ for }x\in\R^{n_1}\mbox{ and }y\in\R^{n_2},\ee
and set $\varphi(x,y)=(g(x,y),y).$ Notice that $g(x,y)=x+\psi(y)(f(x)-x)$, and $f(x)-x$ is bounded, hence
\be \begin{split}g\mbox{ and }\varphi&\mbox{ are Lipschitz, with a bound }L\mbox{ on their Lipschitz norms,}\\
&\mbox{ and } L\mbox{ does not depend on }R.\end{split}\ee

Set $W=\{(x,y)\in\R^{n_1}\times\R^{n_2};f(x,y)\ne (x,y)$, then if $(x,y)\in W$, $\psi(y)\ne 0$ and $f(x)\ne x$. As a result,
\be W\subset B\times B(0,R+1).\ee

We want to compare the measure of $E$ and $\varphi(E)$ in $B\times B(0,R+1)$.

For the part $\varphi(E)\cap (B\times B(0,R))$, observe that by definition of the map $\varphi$, for each $(x,y)\in E,\pi\varphi(x,y)=y$, where $\pi$ is the orthogonal projection on $\R^{n_2}$. This means, $\varphi$ does not change the second coordinate. Hence $\pi[\varphi(E)\cap B\times B(0,R)]\subset E_2\cap B(0,R)$.  

As a result we can apply the coarea formula (c.f \cite{Fe} Corollary 3.3.22) to the set $E$, which is rectifiable because it is minimal.
\be \begin{split}\int_{\varphi(E)\cap (B\times B(0,R))}&||\wedge_{d_2}\pi(x)||dH^d(x)\\
=\int_{E_2\cap B(0,R)}dH^{d_2}(z)&H^{d_1}[\pi^{-1}\{z\}\cap\varphi(E)\cap (B\times B(0,R))].\end{split}\ee

For the left-hand side, notice that $||\wedge_{d_2}\pi|_{\varphi(E)}(x)||=1$ for $x\in\varphi(E)\cap(B\times B(0,R))$, hence
\be \int_{\varphi(E)\cap (B\times B(0,R))}||\wedge_{d_2}\pi(x)||dH^d(x)=H^d(\varphi(E)\cap (B\times B(0,R)).\ee
For the right-hand side, since for any $z\in B(0,R)\cap E_2$, $\pi^{-1}\{z\}\cap\varphi(E)\cap (B\times B(0,R+1))=f(E)\cap B$, therefore
\be H^{d_1}[\pi^{-1}\{z\}\cap\varphi(E)\cap (B\times B(0,R))]=H^{d_1}(E_1\cap B)-c.\ee
As a result,
\be \begin{split}&\int_{E_2}dH^{d_2}(z)H^{d_1}[\pi^{-1}\{z\}\cap\varphi(E)\cap (B\times B(0,R))]\\
=&\int_{E_2\cap B(0,R)}dH^{d_2}(z)H^{d_1}[\pi^{-1}\{z\}\cap\varphi(E)\cap (B\times B(0,R))]\\
=&H^{d_2}(E_2\cap B(0,R))\times [H^{d_1}(E_1\cap B)-c]\\
=&H^d(E\cap(B\times B(0,R))-cH^{d_2}(E_2\cap B(0,R)).
\end{split}\ee

Therefore 
\be H^d(\varphi(E)\cap(B\times B(0,R)))=H^d(E\cap(B\times B(0,R))-cH^{d_2}(E_2\cap B(0,R)).\ee
On the other hand, by definition of $\varphi$, $\varphi(E)\cap [B\times (B(0,R+1)\bs B(0,R))]=\varphi(E\cap [B\times (B(0,R+1)\bs B(0,R))])$, therefore   for any $z\in B(0,R+1)\bs B(0,R)$, by (7.5) 
\be \begin{split}
&H^d(\varphi(E)\cap [B\times (B(0,R+1)\bs B(0,R))])\\
\le& L^dH^d(E\cap [B\times (B(0,R+1)\bs B(0,R))])\\
=& L^dH^{d_1}(E_1\cap B)H^{d_2}(E_2\cap B(0,R+1)\bs B(0,R)),
\end{split}\ee
where $L$ does not depend on $R$.

Combine (7.11) and (7.12), we get
\be \begin{split}H^d(\varphi(E)\cap (B\times B(0,R+1)))\le H^d(E\cap(B\times B(0,R))-cH^{d_2}(E_2\cap B(0,R))\\
+L^dH^{d_1}(E_1\cap B)H^{d_2}(E_2\cap B(0,R+1)\bs B(0,R)).\end{split}\ee
Set $C'=L^dH^{d_1}(E_1\cap B)$, then
\be \begin{split}H^d(E\cap(B\times B(0,R))&-H^d(\varphi(E)\cap (B\times B(0,R+1)))\\
\ge cH^{d_2}(E_2&\cap B(0,R))-C'H^{d_2}(E_2\cap B(0,R+1)\bs B(0,R))\end{split}\ee
where $c,C'$ do not depend on $R$.
Now by the minimality of $E$, we have
\be cH^{d_2}(E_2\cap B(0,R))-C'H^{d_2}(E_2\cap B(0,R+1)\bs B(0,R))\le 0\mbox{ for all }R>0.\ee
So if we can get that 
\be \liminf_{R\to\infty}\frac{H^{d_2}(E_2\cap B(0,R+1)\bs B(0,R))}{H^{d_2}(E_2\cap B(0,R))}=0,\ee
then we can deduce from (7.15) that $c\le 0$ for any $f$ (recall the definition of $c$ in (7.2)), and thus $E_1$ is minimal.
The same argument gives that if
\be \liminf_{R\to\infty}\frac{H^{d_1}(E_1\cap B(0,R+1)\bs B(0,R))}{H^{d_1}(E_1\cap B(0,R))}=0\ee
then $E_2$ is minimal.

We claim that at least one of the (7.16) and (7.17) is true. In fact, denote by $B_i(0,r)=B(0,r)\cap \R^{d_i}$, then for every $R>0$, 
\be \begin{split}(E_2\cap B_2(0,&R+1)\bs B_2(0,R))\times (E_1\cap B_1(0,R+1)\bs B_1(0,R))\\
&\subset E\cap [B_1(0,R+1)\times B_2(0,R+1)]\bs [B_1(0,R)\times B_2(0,R)].\end{split}\ee
But $E$ is minimal, hence by the uniform Ahlfors regularity of minimal sets (cf. \cite{DS00} Proposition 4.1),
\be \begin{split}H^d(E\cap [B_1(0,R+1)\times B_2(0,R+1)]\bs [B_1(0,R)\times B_2(0,R))\\
=o(R^d)=o(H^d(E\cap B_1(0,R)\times B_2(0,R))),\end{split}\ee
which implies that
\be \begin{split}H^d[(E_2\cap B_2(0,R+1)\bs B_2(0,R))\times (E_1\cap B_1(0,R+1)\bs B_1(0,R))]\\
=o(H^d(E\cap B_1(0,R)\times B_2(0,R))),\end{split}\ee
or equivalently,
\be \begin{split}H^{d_1}(E_1\cap B_1(0,R+1)\bs B_1(0,R))\times H^{d_2}((E_2\cap B_2(0,R+1)\bs B_2(0,R))]\\
=o(H^{d_1}(E\cap B_1(0,R))\times H^{d_2}(B_2(0,R))),\end{split}\ee
and hence
\be \lim_{R\to\infty}\frac{H^{d_2}(E_2\cap B(0,R+1)\bs B(0,R))}{H^{d_2}(E_2\cap B(0,R))}\cdot \frac{H^{d_1}(E_1\cap B(0,R+1)\bs B(0,R))}{H^{d_1}(E_1\cap B(0,R))}=0,\ee
which implies that one of the (7.16) and (7.17) is true.

Suppose for example that (7.16) is true, then $E_1$ is minimal. But in this case we can use the Ahlfors regularity of $E_1$ to get (7.17), which implies that $E_2$ is also minimal. Hence both of the two $E_i,i=1,2$ are minimal. \qed

%
% \section{Questions to be answered}
% 
% The proof of the minimality of the cone $Y\times Y$ is part of the classification of singularities for 2-dimensional minimal sets in $\R^n$. Hence immediately after
% 
% Full-length ? existence of a point of $Y\times Y$?  $Y\times T$?

\renewcommand\refname{References}
%\begin{spacing}{1}
\bibliographystyle{plain}
\bibliography{reference}
%\end{spacing}

\end{document}